%% file: LV_bistable_survey.tex
\documentclass[american,a4paper]{amsart}
\pdfoutput=1
\usepackage{refstyle}
\usepackage{mathrsfs}
\usepackage[inline]{enumitem}
\usepackage{mathtools}
\usepackage{amssymb}
\usepackage{graphicx}
\usepackage[unicode=true,pdfusetitle,
 bookmarks=true,bookmarksnumbered=false,bookmarksopen=false,
 breaklinks=false,pdfborder={0 0 1},backref=false,colorlinks=false]{hyperref}
\usepackage{breakurl}
\usepackage{subcaption}
\usepackage{pgfplots}
\usepackage[foot]{amsaddr}
\usepackage{cite}

\newref{fig}{name=Figure~}
\newref{subsec}{name=Subsection~}

\begin{document}

\title[The effect of random dispersal on competitive exclusion]{The effect of random dispersal on competitive exclusion -- A review}

\author{L\'{e}o Girardin}

\thanks{This work was supported by a public grant as part of the Investissement d'avenir project,
reference ANR-11-LABX-0056-LMH, LabEx LMH. This work has been carried out in the framework of 
the NONLOCAL project (ANR-14-CE25-0013) funded by the French National Research Agency (ANR)}
\address{Laboratoire de Math\'{e}matiques d'Orsay, Universit\'{e} Paris
Sud, CNRS, Universit\'{e} Paris-Saclay, 91405 Orsay Cedex, France}
\email{leo.girardin@math.u-psud.fr}

\begin{abstract}
    Does a high dispersal rate provide a competitive advantage when risking competitive exclusion? To this day, the theoretical literature cannot answer this question
    in full generality. The present paper focuses on the simplest mathematical model with two populations differing only in 
    dispersal ability and whose one-dimensional territories are spatially segregated. Although the
    motion of the border between the two territories
    remains elusive
    in general, all cases investigated in the literature concur: either the border does not move at all because
    of some environmental heterogeneity or 
    the fast diffuser chases the slow diffuser. Counterintuitively, it is better to randomly explore the hostile enemy territory, even 
    if it means highly probable death of some individuals, than to ``stay united''. This directly contradicts a celebrated 
    result on the intermediate competition case, emphasizing
    the importance of the competition intensity. Overall, the larger picture remains unclear and the optimal strategy 
    regarding dispersal remains ambiguous. Several open problems worthy of a special attention are raised.
\end{abstract}


\maketitle
\section{Introduction}
The interplay between dispersion and competition is a vast and important problem in theoretical population biology, with applications in
ecology but also in evolution (natural selection precisely originates in the interplay between competitive pressure and 
mutations, namely ``dispersion'' in the phenotypical space), epidemiology (competition between
pathogen strains during spreading epidemics), 
medicine (populations being in this context cell populations). This interplay leads to 
qualitative outcomes (displacement, segregation, etc.) that 
would not appear in the spatially homogeneous, well-mixed counterpart, which makes them
difficult to predict. But such predictions are of the utmost importance, especially since these 
outcomes usually result in some form of spatialized extinction process, which could be a 
goal or on the contrary something to be avoided, depending on the exact biological problem.
An exhaustive overview of the biology and mathematical biology literature on this wide topic is of course impossible; 
the reader is referred for instance to some recent works and references therein
\cite{Cantrell_Cosner_Ruan_2009,Amarasekare_20,Debarre_Lenormand_11,Melbourne_Corn_07,North_Ovaskainen_2007,Phillips2019,Osnas_2015,Berestycki_Zilio_17}.

A common phenomenological mathematical model to study this interplay, inspired by the general population dynamics equation:
\[
    \text{variation in time of the population size} = \text{dispersion} + \text{births} - \text{deaths},
\]
is the deterministic, diffusive and competitive Lotka--Volterra system \cite{Lotka_1924,Volterra_1926,Gause_1934}:
\[
    \begin{cases}
	\partial_{t}N_1 & = \nabla\cdot\left( d_1\nabla N_1 \right) + r_1 N_1\left( 1-\frac{N_1}{K_1} \right)-c_1 N_1 N_2,\\
	\partial_{t}N_2 & = \nabla\cdot\left( d_2\nabla N_2 \right) + r_2 N_2\left( 1-\frac{N_2}{K_2} \right)-c_2 N_1 N_2,
\end{cases}
\]
where $N_1(t,x)\geq 0$ and $N_2(t,x)\geq 0$ are two continuous population densities depending on time $t$ (a real variable) and space $x$ 
(a Euclidean variable $x=\left( x_1,x_2,\dots \right)$), $\partial_t=\frac{\partial}{\partial t}$ is the partial 
derivative with respect to time, $\nabla=\left( \partial_{x_1},\partial_{x_2},\dots \right)$ is the nabla operator (so that $\nabla\bullet$ is the spatial gradient and
$\nabla\cdot\bullet$ is the spatial divergence),
$d_1(t,x)$ and $d_2(t,x)$ are the diffusion (dispersal) rates,
$r_1(t,x)$ and $r_2(t,x)$ are the intrinsic (per capita) growth rates, $K_1(t,x)\geq 0$ and $K_2(t,x)\geq 0$ are the carrying capacities,
$c_1(t,x)\geq 0$ and $c_2(t,x)\geq 0$ are the interpopulation competition rates. 
The dependencies on time and space of the various coefficients account for instance for seasonality or unfavorable regions of space.

Note that in absence of one population, the other grows logistically (no Allee effect) and its density solves the well-known
Fisher--Kolmogorov--Petrovskii--Piskunov partial differential equation (PDE) \cite{Fisher_1937,KPP_1937,Skellam_1951}:
\[
    \partial_{t}N_i = \nabla\cdot\left( d_i\nabla N_i \right) + r_i N_i\left( 1-\frac{N_i}{K_i} \right).
\]
Note also that the paper adopts an ecology vocabulary for the sake of simplicity but, again, the model is abstract and quite general: 
``individuals'' could be cancer cells or infected hosts, etc. 

Assuming for a moment that the environment is spatio-temporally homogeneous (heterogeneous environments will come back later on),
all the coefficients become positive constants and space, time and the two population densities can be nondimensionalized 
to obtain the reduced system:
\begin{equation}
    \begin{cases}
    \partial_{t}u = \Delta u + u\left( 1-u \right)-huv,\\
    \partial_{t}v = d\Delta v + rv\left( 1-v \right)-kuv,
\end{cases}
\label{sys:LV_homo_asym}
\end{equation}
where $d,r,h,k$ are positive constants that can be estimated using field data and where the Fickian diffusion operator reduces
to the simpler spatial Laplacian (corresponding to isotropic Brownian motion of individuals). 
Such a PDE system has indeed been used and discussed
extensively for more than fifty years by modelers in biology. 
A non-exhaustive list illustrating the variety of biological applications includes for instance studies on
the competitive displacement of the red squirrel by the invasive grey squirrel in the British Isles \cite{Okubo_1989}, 
optimization of cancer therapy taking into account the competition between cancer cells that are sensitive to 
the treatment and those that are resistant to it \cite{Carrere_2017,Carrere_2017_2}, 
biodiversity conservation in fire-prone savannas accounting for competition for light and nutrients between 
trees and grass \cite{Yatat_Couteront_Dumont},
reproduction--dispersion trade-offs in experimental bacterial invasions used to study the evolution of 
dispersal \cite{Deforet_2017}.

This synthesis is concerned with results investigating whether the population $v$, in 
order to outcompete the population $u$, 
should have a high or low dispersal rate $d$, all else being equal. 
Here, all else being equal means that 
the two populations only differ in dispersal rate, namely $r=1$ and $h=k$:
\begin{equation}
    \begin{cases}
    \partial_{t}u = \Delta u + u\left( 1-u \right)-kuv,\\
    \partial_{t}v = d\Delta v + v\left( 1-v \right)-kuv.
\end{cases}
\label{sys:LV_homo}
\end{equation}
The symmetry assumption ($h=k$ and $r=1$) prevents pure reaction-driven extinctions that would
strongly perturb the analysis (in other words, cases where, for instance, 
one competitor feels much less pressure than the other and might prevail despite a poorly chosen
dispersal strategy are discarded).
This assumption will simplify a lot the forthcoming presentation, although many results actually 
remain true under specific, yet more general, assumptions on $h$ and $r$.

In spatially homogeneous, well-mixed cases, the system (\ref{sys:LV_homo}) is strongly determined
by the sign of $k-1$ (\textit{e.g.}, \cite[Chapter 7, Section 7.9]{Iannelli_Pugliese}). 
On one hand, in the weak competition case ($k<1$), the system
is systematically driven to coexistence. On the contrary, in the strong competition case ($k>1$), 
both populations are able to wipe out the other provided they are numerically sufficiently superior. 
The intermediate case ($k=1$), corresponding to a competitive pressure exerted completely blindly, is degenerate 
(all pairs $(u,v)$ satisfying $u+v=1$ are steady states) and is usually discarded. Note that ``blind'' competition means here
that one individual competes uniformly with all surrounding individuals, independently of the population to 
which they belong; the population label of competitors is in some sense not seen, not taken into account.

In spatially structured systems, the picture is more complicated.
A very important paper of mathematical biology, due to Dockery, Hutson, Mischaikow and Pernarowski 
\cite{Dockery_1998}, established that the blind competition case $k=1$ becomes relevant in spatially
heterogeneous environments with an intrinsic growth rate $a(x)$:
\[
    \begin{cases}
    \partial_{t}u = \Delta u + u\left( a -u \right)-uv,\\
    \partial_{t}v = d\Delta v + v\left( a -v \right)-uv.
\end{cases}
\]
Assuming that the domain where the populations evolve is bounded with no-flux boundary conditions (say, an island or a Petri dish), 
the authors showed that $v$ wipes out $u$ whenever it is the slower diffuser, $d<1$. Of course, by symmetry, $u$ wipes out $v$ if $d>1$. 
This was interpreted as follows: because of the interplay between heterogeneity and competition, it is a better strategy to 
claim favorable areas and to defend them collectively by remaining concentrated there than to randomly explore unfavorable
areas, where deaths due to the environment are more likely. In other words, $v$ wins if and only if, compared to $u$, 
its individuals remain ``united'' instead of venturing alone in unknown areas. In the present paper,
such a result is referred to as a ``Unity is strength''-type result. 
The analysis of Dockery \textit{et al.} relied entirely upon the monostability of the system induced by the combination of
spatial heterogeneity and difference in diffusion rates: the only stable steady state is the one where the slow diffuser
persists while the fast diffuser vanishes. Initial conditions do not matter: even if, 
initially, the fast diffuser is vastly superior numerically, the slow diffuser will eventually 
prevail. 

But what if both semi-extinct steady states are stable, so that the stability analysis does not
suffice to conclude and initial conditions matter? 
As explained above, bistability is for instance achieved in spatially 
homogeneous systems (\ref{sys:LV_homo}) with strong competition ($k>1$). 
With such systems, the success of a dispersal strategy is a more
delicate notion that can be defined in a few ways. 

For instance, the diffusion-induced extinction property could be used to define
this success. This criterion uses initial conditions that are exactly calibrated so that neither $u$ nor $v$ takes over in the
absence of diffusion or with equal diffusion rates ($d=1$). 
Given such initial conditions, what is the outcome when taking the unequal diffusion into account? 
In homogeneous environments, where this balance condition simply reads $u(0,x)=v(0,x)$ at every $x$,
and provided the habitat is one-dimensional with no-flux boundary conditions and 
the interpopulation competition rate $k$ is equal to $2$, Ninomiya \cite{Ninomiya_1995} 
showed that there exist values of $d$ larger than $1$ but close to it and carefully chosen 
initial conditions satisfying the above condition such that $v$ wipes out $u$. 
The fast diffuser prevails: in this sense, ``Unity is not strength'' (one could even say ``Disunity is strength''). 

Nevertheless, this definition of success is unsatisfying, as it uses very precise initial conditions that are in 
some sense artificial and would not appear in the nature. Is there a more robust and natural definition (that might 
\textit{a priori} disagree with the conclusions of Ninomiya \cite{Ninomiya_1995} and agree with those of Dockery 
\textit{et al.} \cite{Dockery_1998})?

The strong competition case $k>1$ is also known in the mathematical
ecology literature as the competitive exclusion case \cite{Gause_1934}: persistence of both species 
can occur only if the two niches are differentiated. 
In the setting of this paper, niches are purely geographical, and their
differentiation means that, roughly speaking, $u$ is positive where $v$ is close to $0$ and
vice-versa (note the sharp contradistinction with Ninomiya's balance condition $u(0,x)=v(0,x)$).
If the territories are segregated, then borders between these territories naturally arise.
At these interfaces, the two populations meet frequently and individuals
compete fiercely to chase competitors and take over. In this context, it seems natural
to track the motion of the interface and to define a dispersal strategy as successfull if it 
leads to taking over the territory of the opponent, namely to territorial expansion. 

This definition agrees with situations studied in the biological literature 
(\textit{e.g.}, \cite{Mitani_Watts_Amsler,Brown_71,Bohn_Amundsen_,Berestycki_Zilio_17}).

\begin{figure}
        \resizebox{.8\linewidth}{!}{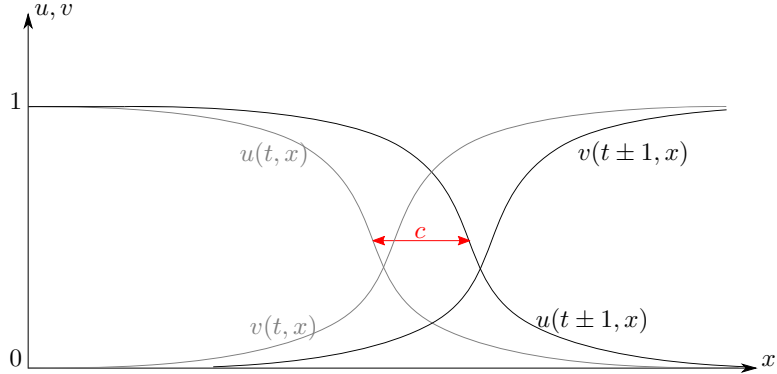}
	\caption{Bistable traveling wave solution of the system (\ref{sys:LV_homo}) with speed $c$. Depending on the sign of $c$,
    the black curves correspond to the time $t+1$ or $t-1$.}
    \label{fig:traveling_wave}
\end{figure}

Mathematically, this definition translates in homogeneous environments into the study of a particular solution of the 
system (\ref{sys:LV_homo}), referred to as
a traveling wave, that has a constant profile and a constant speed and evolves in the infinite real line (approximating a very large
one-dimensional domain where propagation phenomena matter). This solution is illustrated in \figref{traveling_wave}. Its existence
and its uniqueness were confirmed in the '80s and '90s \cite{Gardner_1982,Kan_on_1995,V_V_V}. In this context,
the success of the dispersal strategy is simply given by the sign of the speed of the wave. However, in contrast with
the existence and uniqueness of the wave, this sign is in general a
very difficult mathematical problem, that cannot be solved by any standard tool of the analysis of PDEs. Only partial results are
known and these are the main topic of this synthesis paper. 
It turns out that they are all in agreement with Ninomiya \cite{Ninomiya_1995}: 
in situations of competitive exclusion due to strong interpopulation competition, ``Unity is not strength''.

The paper is organized as follows.
Section 2 gives an exhaustive survey of these partial results, the last subsection being
devoted to a delicate extension in spatially periodic media. 
Section 3 lists some open problems that should, in my opinion, attract the attention of the community.

\section{Known results}

In this section the known results are presented.

\subsection{Homogeneous environment}

In this subsection the focus is on homogeneous environments (constant coefficients), so that the system
is (\ref{sys:LV_homo})
and the interesting solution is indeed the traveling wave solution of \figref{traveling_wave}. Its speed is denoted
$c_{k,d}$. By symmetry and without loss of generality, the assumption $d>1$ stands, so that ``Unity is strength'' holds true
if and only if $c_{k,d}>0$.

Note that the opposite invasion ($u$ on the right, $v$ on the left, with speed $\tilde{c}_{k,d}=-c_{k,d}$) could be considered 
and is in fact the one considered in some of the forthcoming references.
Without loss of generality, the only case considered hereafter is that of \figref{traveling_wave}.

\begin{figure}
        \resizebox{.9\linewidth}{!}{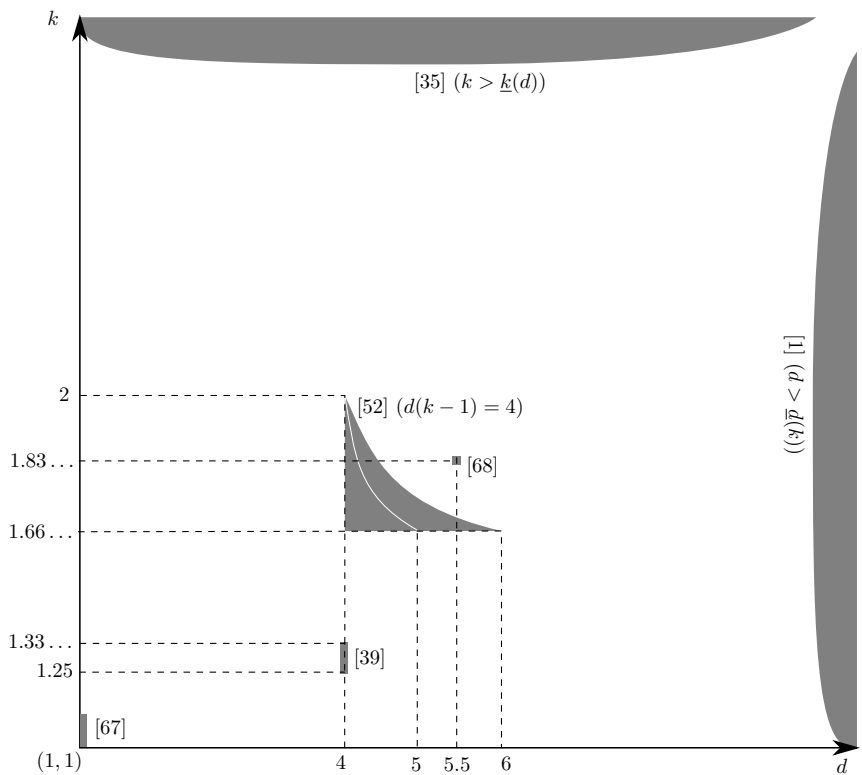}
	\caption{Localization in the $(d,k)$ plane of the known results for the sign of the 
	    traveling wave speed. 
	    Shaded areas: $c_{k,d}<0$; white areas: unknown sign; on the $k$ axis: $c_{k,1}=0$. 
	    The (in)equations describe the non-straight lines.}
    \label{fig:homogeneous_env}
\end{figure}

All these results are summarized in \figref{homogeneous_env}.

\subsubsection{Very special choices of parameters}

In 2001, Rodrigo and Mimura \cite{Rodrigo_Mimura_2001} computed nine exact families of traveling wave 
solutions of the system (\ref{sys:LV_homo_asym}) by looking for closed-form solutions 
(a hyperbolic tangent ansatz) and compatibility conditions on the parameters $d,r,h,k$. 
A comparison of their results with the symmetry assumption ($r=1$, $h=k$) yields
\[
    c_{11/6,11/2}=-\frac{\sqrt{6}}{12}<0.
\]

In 2013, using the strict monotonicity of the wave speed of the system (\ref{sys:LV_homo_asym}) with respect to 
$h$ and $k$, established in 1995 by Kan-on \cite{Kan_on_1995}, and an analysis
of the special case where the wave speed is zero, Guo and Lin \cite{Guo_Lin_2013} stated a few algebraic conditions on the 
parameters $d,r,h,k$ sufficient to characterize the sign of the speed. A comparison of their results with the symmetry assumption yields
\[
    c_{k,4}<0\quad\text{if }\frac{5}{4}\leq k\leq \frac{4}{3}.
\]

In 2019, using comparison arguments and super- and sub-solutions, Ma, Huang and Ou \cite{Ma_Huang_Ou_2019} obtained five new algebraic
conditions on the parameters $d,r,h,k$ sufficient to characterize the sign of the speed. With the symmetry assumption, their conditions
read
\[
    c_{k,d}<0\quad\text{if }k>\frac{5}{3}\text{ and }d\in\left( 4,\frac{2k}{k-1} \right)\cup\left( \frac{2k}{k-1},\frac{4}{k-1} \right).
\]
Note that the case $d=\frac{2k}{k-1}$ is not solved (see \figref{homogeneous_env}).

\subsubsection{Almost equal dispersal rates and almost blind competition}

In 2017, Risler \cite{Risler_2017} investigated approximations of $c_{1+(\delta k)^2,1+\delta d}$.
It turns out that there is a singularity at $k=1$. When $k>1$, the wave speed is a regular function of $(k,d)$
and, using the identity $c_{k,1}=0$, the following third-order (in $\delta d$ and $\delta k$) approximation holds true
(all partial derivatives being implicitly evaluated at $(k,1)$ for clarity):
\[
    c_{k+(\delta k)^2,1+\delta d} \sim \left[\frac{\partial c}{\partial d}
	+ \frac{1}{2}\frac{\partial^2 c}{\partial d^2}\delta d
	+ \frac{1}{6}\frac{\partial^3 c}{\partial d^3}(\delta d)^2
        + \frac{\partial^2 c}{\partial k\partial d}(\delta k)^2\right]\delta d.
\]
By a singular perturbation approach, Risler managed to prove that the quantity
$\frac{1}{\delta k}\left(\frac{\partial c}{\partial d} + \frac{\partial^2 c}{\partial k\partial d}(\delta k)^2\right)$
converges, as $k\to 1^{+}$, to some negative constant $A<0$.
As a consequence, in the parameter regime $0<\delta d\ll\delta k\ll 1$ where the second and third order terms in $\delta d$
can be safely ignored,
\[
    c_{1+(\delta k)^2,1+\delta d}\sim -|A|\delta k\delta d<0.
\]

This result is especially interesting since, compared with the result of Dockery \textit{et al.} \cite{Dockery_1998},
it shows a trade-off at $d=1+\delta d$, $k=1$ between the intensity of the interpopulation competition and the heterogeneity 
of spatial resources.
On one hand, changing the homogeneous resources into spatially heterogeneous resources $a(x)=1+\delta a(x)$ favors the 
slow diffuser $u$; but on the other hand, adding a small bump in interpopulation competition $(\delta k)^2$ favors the fast diffuser $v$.

\subsubsection{Large gap between the dispersal rates}

In 2004 and 2005, Heinze \textit{et al.} \cite{Heinze_et_al_2004,Heinze_Schweiz} studied the degenerate system (\ref{sys:LV_homo_asym})
where one of the two populations does not diffuse at all ($d=0$). 
Inspired by this work, Alzahrani, Davidson and Dodds \cite{Alzahrani_Davidson_Dodds_2010} studied in 2010
the strong dispersal limit $d\to +\infty$ of the system (\ref{sys:LV_homo_asym}) (which is equivalent,
up to a change of variable, to the limit $d\to 0$). 
In this regime and taking into account the symmetry assumption, they proved 
by energy methods that $c_{k,d}/\sqrt{d}$ converges to a finite limit $l_{k,\infty}\in\left[ -2,0 \right)$.
Consequently, for any $k>1$, there exists $\underline{d}(k)\geq 1$ such that
\[
    c_{k,d}<0\quad\text{if }d>\underline{d}(k).
\]
The variations (monotonicity, convexity) and limits as $k\to 1$ or $k\to+\infty$ of the optimal function $\underline{d}$ are unknown.
The graph of the function $\underline{d}$ in \figref{homogeneous_env} is an arbitrary choice.

\begin{figure}
        \resizebox{.9\linewidth}{!}{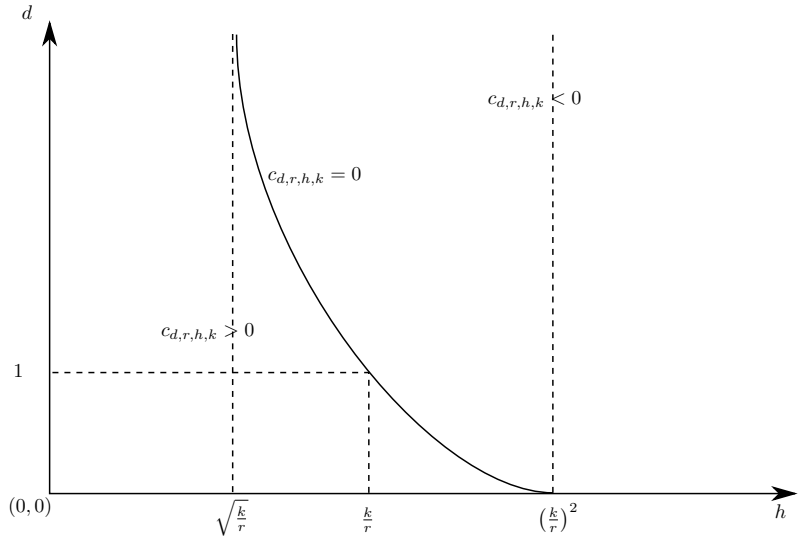}
	\caption{Illustration of a global ``Disunity is strength''-type conjecture 
	    by Alzahrani, Davidson and Dodds 
	    \cite{Alzahrani_Davidson_Dodds_2012}. Here, the parameters $r$ and $k$ of
	    the system (\ref{sys:LV_homo_asym}) are
	    arbitrarily fixed and the traveling wave speed is studied as a function of $(h,d)$.
	    The fact that it is zero at the point $\left(\frac{k}{r},1\right)$ is obvious; the
	    fact that it is asymptotically zero at
	    $\left(\sqrt{k/r},+\infty\right)$ and $\left(\left(k/r\right)^2,0\right)$ 
	    is proved in \cite{Alzahrani_Davidson_Dodds_2010}. 
	    The rest, and in particular the fact that the $0$-level set is a monotone graph, 
	is only supported by numerical experiments.}
    \label{fig:ADD_2012}
\end{figure}

In 2012, the same authors published a sequel \cite{Alzahrani_Davidson_Dodds_2012}
where they numerically completed the picture for intermediate values of $d$ and stated
a global ``Disunity is strength''-type conjecture. This conjecture is summarized 
in \figref{ADD_2012}, which is
basically an adaptation of \cite[Figure 6]{Alzahrani_Davidson_Dodds_2012} to the present setting and 
notations.

\subsubsection{Very strong interpopulation competition}

In 2015, Girardin and Nadin \cite{Girardin_Nadin_2015} studied the very strong competition limit $k\to+\infty$. In this regime,
the two territories are completely segregated ($uv=0$ everywhere).
Using this property, the authors managed to prove that $c_{k,d}$ converges to a limit $c_{\infty,d}\in\left( -2\sqrt{d},0 \right)$. 
Therefore, for any $d>1$, there exists $\underline{k}(d)\geq 1$ such that
\[
    c_{k,d}<0\quad\text{if }k>\underline{k}(d).
\]
Note that the identity $c_{k,1}=0$ leads to $\lim_{d\to 1}\underline{k}(d)=+\infty$.
Apart from this limit, the graph of the optimal function $\underline{k}$ is unknown and the choice in \figref{homogeneous_env} is 
arbitrary.

\subsection{Spatially periodic environment}

Recalling the seminal paper of Dockery \textit{et al.} \cite{Dockery_1998}, it is quite natural to investigate a possible extension
of the preceding ``Disunity is strength''-type results to spatially
heterogeneous environments. In heterogeneous media, traveling waves as defined above (eternal solutions with constant speed and constant
profile) are not solutions anymore, yet generalizations of the concept 
exist (for instance, the generalized transition fronts \cite{Berestycki_Ham-1,Berestycki_Ham-2}). 
Even so, the aforementioned extension turns out to be an impossible task in 
full generality: heterogeneities can block bistable invasions (a phenomenenon referred to as pinning, quenching or blocking, where 
the interface between the territories stops moving and both niches persist) 
\cite{Barton_Turelli,Matano_79,Matano_Mimura_83,Berestycki_Bouhours_Chapuisat} and even
repel them \cite{Zlatos_2017}. 
Therefore assumptions have to be made prior to any meaningful extension.

To the best of my knowledge, the only special form of heterogeneity whose effect on the motion
of the strongly competitive interface has been investigated is the spatially periodic one.
It is natural to consider such a form of heterogeneity for mainly two reasons. 
First, it has convenient mathematical properties, thanks to 
which it is reasonable to expect a clear result. Second, it is heterogeneous 
everywhere -- contrarily, say, to an homogeneous domain with one obstacle where individuals 
cannot go or survive -- and therefore the effect of the heterogeneity on the motion 
of the interface is never attenuated due to an increasing distance. In this sense, the 
analysis of the periodic framework is an important step forward. Accordingly,
researchers have been interested in the invasion of a single population in a periodic 
environment for several decades 
(\textit{e.g.}, \cite{Gartner_Freidlin,Kawasaki_Shige,Xin_1993} and more recently 
\cite{Berestycki_Ham_1,Berestycki_Ham_2,Berestycki_Ham_3}), however studies on competitive
displacements started much more recently. Again, the stability analysis is the starting point;
since the identification of the invader remains quite direct in monostable cases 
(on which the reader is referred for instance to \cite{Yu_Zhao_2015,Maciel_Lutscher_18} and references therein), 
it is again natural to focus on cases with two -- or more, see below -- stable steady states.

In the periodic framework, traveling waves are replaced by solutions known as traveling pulsating (or periodic) waves (or fronts) 
\cite{Berestycki_Ham_3,Kawasaki_Shige} (pulsating waves hereafter). Such a solution still has a
constant speed $c_{k,d}$ but its profile in the moving frame with speed $c_{k,d}$ varies slightly
as time goes on (illustrations can be found for instance in \cite{Kawasaki_Shige}).
In contradistinction with the homogeneous case, the existence of pulsating waves is not systematic:
stable periodic coexistence steady states might now exist and induce repelling phenomena preventing 
the formation of wave patterns. 
A very recent general result by Du, Li and Wu \cite{Du_Li_Wu_2018} provides
conditions on some eigenvalues ensuring the existence of bistable pulsating waves in one-dimensional
spatially periodic media, however it remains
difficult to convert these abstract conditions into practically verifiable conditions.

\subsubsection{Oscillating diffusion rate}
In a recent paper, Hutridurga and Venkatamaran \cite{Hutridurga_Ven}
considered a system where the reaction term is homogeneous but where the diffusion rate $d(x)$ is 
spatially periodic in a one- or two-dimensional environment:
\begin{equation*}
    \begin{cases}
	\partial_{t}u = \Delta u + u\left( 1-u \right)-kuv,\\
	\partial_{t}v = \nabla\cdot\left( d\nabla v \right) + v\left( 1-v \right)-\alpha kuv.
\end{cases}
\end{equation*}

They fixed $\alpha<1$, so that if $d$ is identically equal to $1$ then $v$ has a competitive 
advantage and chases $u$. Fixing then $d$ of the form 
\[
    d_k(x)=1+\frac{3}{4}\sin\left( 2k\pi x \right)
\]
and considering the regime $k\gg 1$ 
(thus both the interpopulation competition rate and the frequency of the oscillations are very 
large), they numerically observed invasion 
reversals: a uniform diffusion rate seems to confer a competitive advantage over an oscillating one 
with the same mean value. Provided this remains true if the mean value of $d$ is slightly 
increased (say, by continuity of $c_{k,d}$), then it yields
a noticeable case of ``Unity is strength''-type result in the strongly competing regime.

\subsubsection{Oscillating reaction terms}
Girardin \cite{Girardin_2016}, together with Nadin \cite{Girardin_Nadin_2016} and Zilio 
\cite{Girardin_Zilio}, published from 2017 to 2019 a series of three articles investigating
analytically the very strong competition limit $k\to+\infty$ of a system with spatially 
homogeneous diffusion rates but with quite general spatially periodic reaction terms in 
one-dimensional environments. 
In particular, the form of spatially heterogeneous reaction terms used by Dockery and his 
collaborators \cite{Dockery_1998} is included in the setting of the 
first and second parts of the series \cite{Girardin_2016,Girardin_Nadin_2016} but, for 
technical reasons, not in the setting of the third part \cite{Girardin_Zilio}.
The presentation below is simplified by assuming 
that the system is now, for some uniformly positive periodic distribution of resources $\mu(x)$,
\begin{equation}
    \begin{cases}
	\partial_{t}u = \partial_{x}^2 u + \mu\left[ u\left( 1-u \right)-kuv\right],\\
	\partial_{t}v = d\partial_{x}^2 v + \mu\left[ v\left( 1-v \right)-kuv\right].
\end{cases}
\label{sys:LV_homo_periodic}
\end{equation}
Nevertheless, keep in mind that many results of the series remain true with more general systems.

In the first part of the series \cite{Girardin_2016}, it was proved that a bistable pulsating 
wave exists provided the interpopulation competition rate $k$ and the frequency of the environment 
$\mu$ are both sufficiently large. This existence result, consistent with the aforementioned abstract
one \cite{Du_Li_Wu_2018}, confirms the general principle according
to which higher dispersal distances destabilize coexistence \cite{Boeye_2014,Debarre_Lenormand_11}.

In the second part \cite{Girardin_Nadin_2016}, it was proved that, provided a pulsating wave 
exists, then in the very strong competition limit $k\to+\infty$, a ``Disunity is strength''-type 
result holds true: there exists a threshold $d^\star\geq 1$ such that, if $d>d^\star$, 
$c_{\infty,d}<0$, whence in particular 
$c_{k,d}<0$ if $k$ is large enough. Whatever the exact shape of $\mu$ is, the fast diffuser $v$ still has a competitive advantage. 

But where does this threshold $d^\star$ come from? What is its biological meaning? When $1<d\leq d^\star$, the limit $c_{\infty,d}$ 
is nonpositive and can actually be zero. This does not mean that $c_{k,d}$ is zero for large values of $k$,
but it remains unclear whether the convergence is from below ($c_{k,d}<0$ if $k\gg 1$), from above ($c_{k,d}>0$ if $k\gg 1$), 
stationary ($c_{k,d}=0$ if $k\gg 1$) or oscillating (changing signs). Assuming for a moment that $c_{k,d}=0$ is realized in some cases,
then these are cases of coexistence by segregation in a wave-like pattern (pinning phenomenon) with a large but finite (thus more realistic)
interpopulation competition rate.

\begin{figure}
        \resizebox{.8\linewidth}{!}{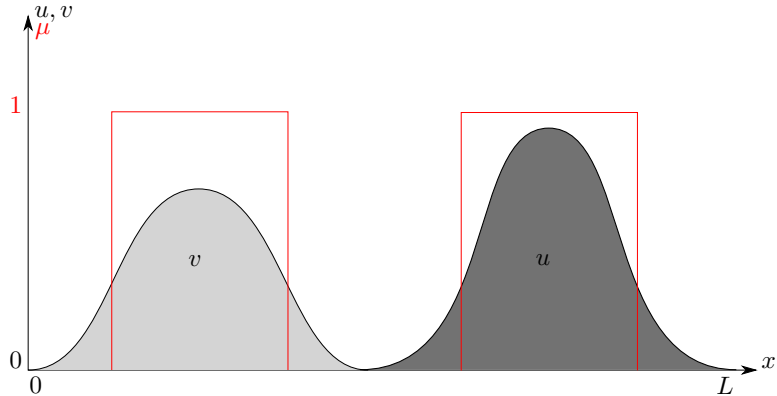}
	\caption{A stable periodic segregated stationary state as constructed by Girardin
	    and Zilio \cite{Girardin_Zilio}. Here the periodic functions $\mu,u,v$ are only 
	    presented in one spatial period $[0,L]$.}
	    \label{fig:GZ}
\end{figure}

In the third part \cite{Girardin_Zilio}, the possibility of wave repelling was explored further. The authors considered 
the specific case of a patchy environment, in which favorable patches where resources are 
abondant and homogeneous ($\mu\simeq 1$) are separated by 
neutral patches where only dispersion occurs ($\mu\simeq 0$).
Assuming a very large (but finite) competition rate $k$ and a small enough frequency of the environment (\textit{i.e.} large patches
or small dispersal distances), they proved the existence of a stable periodic segregated stationary state describing a situation
where the population $u$ settles in the oddly numbered favorable patches while the population $v$ settles in the evenly numbered
favorable patches (see \figref{GZ}). Although such a steady state solution of the system (\ref{sys:LV_homo_periodic}) is unable to induce
a repelling phenomenon (\textit{i.e.} a pulsating wave still exists, at least in the limit $k\to+\infty$ \cite{Nolen_Ryzhik_09}),
the exact same result applies if resources are \textit{a priori} specialized (the oddly numbered, respectively evenly numbered, 
favorable patches are favorable only for $u$, respectively $v$) and for such resources, the existence of pulsating waves 
remains an open problem.

Put together, the second and third part of the series show how spatially periodic reaction 
terms influence competitive displacements when the interpopulation competition is very strong:
the invasion of the fast diffuser $v$ into the territory of $u$ might be blocked or repelled by a
coexistence state but it will not be reversed 
(up to the likely uncommon cases where $1<d\leq d^\star$ and $0<c_{k,d}\ll 1$, mentioned earlier).

\section{Open problems}
In this section, directions of research worthy of a special attention from the community are raised. Although the focus remains on 
theoretical questions and viewpoints, comparisons with \textit{in vitro} models and field
data are also obviously interesting and important.

\subsection{Complete picture in homogeneous environments}
Of course, in view of \figref{homogeneous_env} and \figref{ADD_2012}, it is tempting to 
conjecture that $c_{k,d}<0$ is globally true as soon as $d>1$.
As such a result would be in sharp contradistinction with the widely accepted 
``Unity is strength''-type conclusions of Dockery 
\textit{et al.} \cite{Dockery_1998} and would therefore confirm the trade-off between 
interpopulation competition and spatial 
heterogeneity first hinted by Risler \cite{Risler_2017}, this is definitely the most 
important open problem raised by the present paper.

In order to support further this conjecture, a series of numerical experiments was conducted during the writing of the
present paper \footnote{Simulations were run in \textit{GNU Octave} \cite{Octave}. Details are given in the appendix.}.
These are presented in \figref{numerics}. In view of these experiments, the conjecture seems to be correct. 
In addition, a form of monotonicity with respect to both $k$ and $d$ appears in \figref{numerics}. The 
monotonicity with respect to $d$ was already conjectured and discussed at length in 2012 by
Alzahrani, Davidson and Dodds \cite{Alzahrani_Davidson_Dodds_2012}. The monotonicity with 
respect to $k$ is, to my knowledge, a new conjecture.

\begin{figure}
    \begin{subfigure}{.95\linewidth}
	\resizebox{\linewidth}{!}{\input{heatmap.tex}}
	\caption{Heat map}
    \end{subfigure}

    \begin{subfigure}{.95\linewidth}
	\resizebox{\linewidth}{!}{\input{level_sets.tex}}
    	\caption{$0,-0.1,-0.2,\dots,-1.2$-level sets (the $0$-level set, in white, is a dashed line for clarity; 
	the lack of smoothness is a numerical artifact)}
    \end{subfigure}
    \caption{Numerical approximations of the traveling wave speed $c_{k,d}$, presented in the $(d,k)$ plane and
    color-coded according to the bar on the right-hand side of panel (A)}
    \label{fig:numerics}
\end{figure}
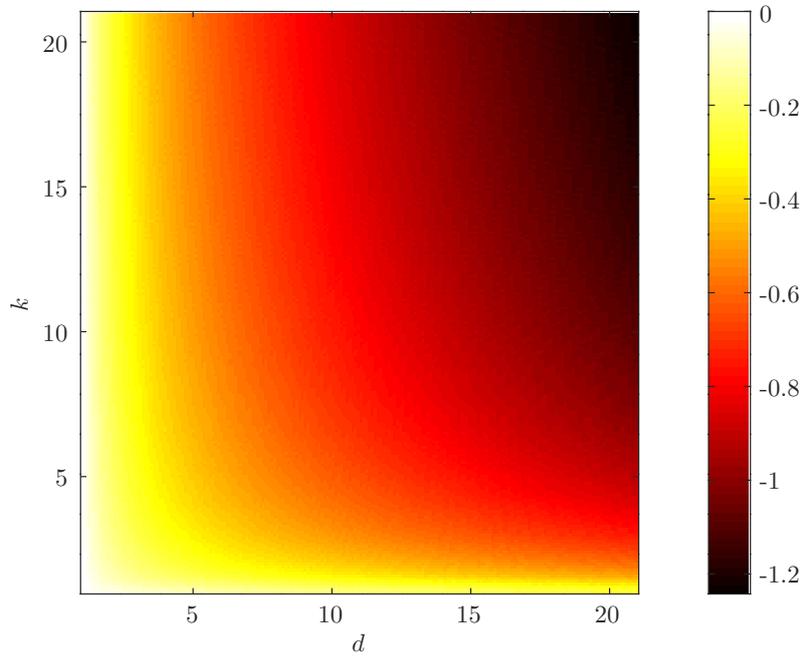
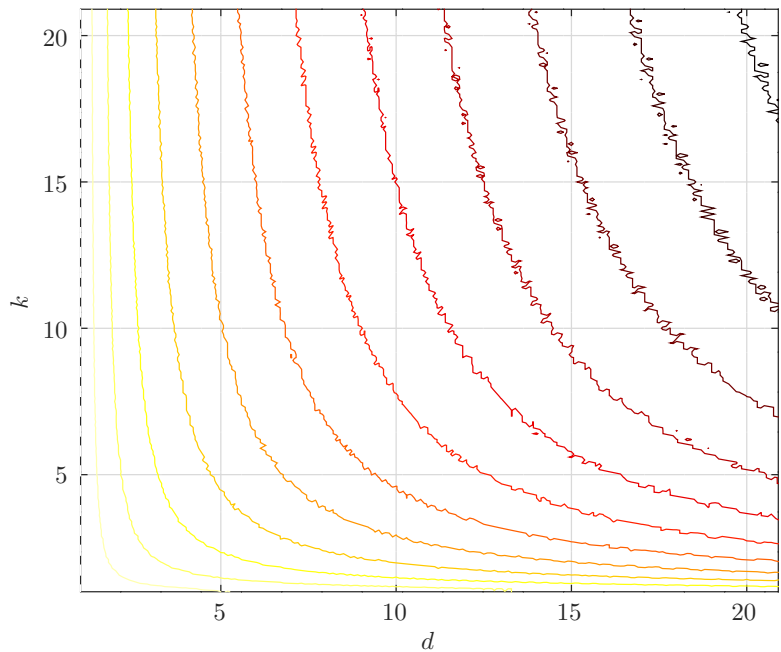

A rigorous analytical proof of the global ``Disunity is strength'' theorem remains to this day a completely open and 
challenging mathematical problem. 

\subsection{Better understanding of the coexistence by segregation in spatially periodic environments}
As explained above, the cases where $d>1$ and $c_{\infty,d}=0$ in spatially periodic media should definitely be investigated further. 
How does the convergence as $k\to+\infty$ occur? In particular, is it possible to observe a pinning phenomenon, 
namely a wave-like coexistence steady state, with biologically relevant finite values of $k$? 
Of course, this leads to the very natural question of whether, in spatially periodic media, ``Unity is strength'' can be achieved
away from the strong competition regime, say close to $k=1$. To the best of my knowledge, this is a completely open problem.

Note that the literature on scalar bistable reaction--diffusion equations tends to indicate that direct constuctions might directly
yield the following results in the very strongly competing regime:
\begin{itemize}
    \item a pinning phenomenon for a system with two oscillating diffusion rates and two uniform reaction terms \cite{Xin_1993};
    \item a repelling phenomenon for a system with two uniform diffusion rates and two oscillating reaction terms \cite{Zlatos_2017}.
\end{itemize}

\subsection{Extension to more general heterogeneous environments}
Extending the spatially periodic ``Unity is not strength''-type results of Girardin \textit{et al.} 
\cite{Girardin_2016,Girardin_Nadin_2016,Girardin_Zilio} in two- or three-dimensional environments (periodic tilings) 
should be possible, to a certain degree at least, but difficult mathematical obstacles will arise.

The system (\ref{sys:LV_homo_periodic}) still makes sense if the amount of resources $\mu$ is sign-changing, namely if some areas
are so unfavorable that they are actually deadly even in absence of competition. Periodic environments with deadly areas
are more subtle and form a very interesting problem. Indeed, at least intuitively, the more lethal the unfavorable areas are, 
the more likely ``Unity is strength'' seems. Applying formally the formulas of the uniformly favorable case 
\cite{Girardin_Nadin_2016}, a condition for $c_{\infty,d}>0$ arises: the negativity of the mean value of $\mu$. Although such
an environment is obviously unfavorable and even harmul, it does not automatically mean that neither $u$ nor $v$ can settle in it.
In fact, it is known from earlier studies that if $\mu$ has a specific form, then both $u$ and $v$ will be able to settle in 
absence of the competitor despite the negativity of the mean value of $\mu$ \cite[Theorem 2.12]{Berestycki_Ham_1}.
Hence the negativity of the mean value of $\mu$ is not a contradictory condition in itself: a pulsating wave describing the
competition between $u$ and $v$ might still exist. In other words, it seems that a ``Unity is strength''-type result
is within reach: such a problem needs to be investigated more thoroughly (first, existence of the wave; second,
sign of the wave speed; the pre-existing methods \cite{Girardin_2016,Girardin_Nadin_2016} will not be sufficient as
they require the positivity of the minimum of $\mu$ but with care it might be possible to build from them). 

In temporally periodic environments, the existence of a unique bistable pulsating wave was proved
in 2013 by Bao and Wang \cite{Bao_Wang_2013} under some conditions on eigenvalues (reminiscent of 
those for spatially periodic environments \cite{Du_Li_Wu_2018}). 
However the sign of the wave speed is completely unknown and, in particular, it is unclear whether
the approach developed by Girardin and Nadin to study the very strong competition limit $k\to+\infty$ 
\cite{Girardin_Nadin_2015,Girardin_Nadin_2016}
can be used again (in particular, the sign of the wave speed cannot be obtained anymore by integration by 
parts of the limiting equation \cite[Proposition 3.7]{Girardin_Nadin_2016}). 
\textit{A fortiori}, spatio-temporally periodic environments are even more elusive. 
Note that it is known that the introduction of temporal periodicity in the 
spatially heterogeneous model of Dockery \textit{et al.} \cite{Dockery_1998} provides in some specific cases a 
competitive advantage to the fast diffuser (``Unity is not strength''-type result) and thus compensates
the effect of the spatial heterogeneity \cite{Hutson_Mischaikow_Polacik}.

Apart from the periodic regime, the effect of environmental heterogeneities on the motion of
the interface is an entirely open problem. In view of the reaction--diffusion literature,
other special regimes (almost periodic, ergodic, finite number of obstacles, cylindrical domains, etc.) can reasonably be considered. 
However this requires a much more involved mathematical analysis, with very specified techniques depending strongly on
the choice of special regime.
In particular, regarding spatial domains with boundaries, the existence of 
a stable steady state of the system with segregated niches strongly depends on 
the exact shape of the domain 
(\textit{e.g.}, \cite{Kishimoto_Wein,Matano_Mimura_83,Kan-on_Yanagida_1993}).

\subsection{Comparison between the Brownian motion and other motion strategies}
Finally, it is important to assert the dependency of the aforementioned results on the exact dispersal strategy of the populations. 
The use of diffusion to model animal movement in ecology, for instance, is debated by 
a wide literature and alternative models exist (\textit{e.g.}
\cite{Benichou_2011,Giuggioli_2010,Ovaskainen_2008,Cantrell_Cosner_Ruan_2009} and references therein). 

More specifically, in what follows we briefly discuss three possibilities: density-dependent dispersal, 
resources-dependent dispersal and long-range dispersal. Note that density-dependent dispersal operators 
can be coupled (when the dispersal strategy depends on the density of the competitor) or not.

Potts and Petrovskii \cite{Potts_Petrovskii} recently studied numerically the 
influence of a density-dependent coupled operator with an aggressive taxis term pushing
the slow diffuser $u$ towards the fast diffuser $v$. Interestingly, they found a special range of parameters in which
$u$ chases $v$, thus reversing the conclusion without taxis. The authors refer to their result as a ``Fortune favors the bold''-type
result: boldness (aggressive advection) can compensate unity (low diffusion rate). An analytical proof of such a result would be
interesting. Similarly, the effect of cross-diffusion and self-diffusion (\textit{e.g.}, \cite{SKT_79,Matano_Mimura_83,Dreher_2007}) could
be investigated, although mathematically it is a very challenging problem (to my knowledge, even the existence of traveling waves is unclear) 
and at the moment it seems that only perturbative results or numerical experiments are within reach. Intuitively, following the 
``Disunity is strength'' logic, it could be conjectured that, opposed to standard diffusion, self-diffusion is a losing strategy 
whereas cross-diffusion is a winning one.

Regarding resources-dependent dispersal, many models could be considered; let me suggest for instance conditional dispersal and ideal free 
dispersal \cite[Chapter 11]{Cantrell_Cosner_Ruan_2009}.

Finally, in ecological or evolutionary contexts, it is very natural to replace the Laplacian diffusion operator by a nonlocal 
operator accounting for rare long-range dispersal events (\textit{e.g.}, \cite{Benichou_2011,Kimura_1965,Hutson_Martinez_2003}).
Intuitively, such dispersion operators might be roughly understood
as diffusion operators with huge diffusion rates \cite{Hutson_Martinez_2003}. Testing such 
an intuition against the competitive model is therefore very tempting, and
in particular it would be natural to study the case where $v$ still diffuses locally with a rate $d$ but $u$ diffuses nonlocally with a 
normalized rate: in such a case, does $v$ always lose, whatever the value of $d$?

\subsection{An evolutionary trap?}
The ``Disunity is strength''-type results presented above and the ``Unity is strength''-type evolutionary result
can be confronted in an interesting way: what if a population subjected to evolution of dispersal encountered a strong competitor?
One possible mathematical model for this problem is the following system:
\begin{equation*}
    \begin{cases}
    \partial_{t}u = \Delta u + u\left( a-u \right)-hu\left( v_1+v_2 \right),\\
    \partial_{t}v_1 = d_1\Delta v_1 + v_1\left( a-v_1-v_2 \right)-kuv_1 + \alpha(v_2-v_1),\\
    \partial_{t}v_2 = d_2\Delta v_2 + v_2\left( a-v_2-v_2 \right)-kuv_1 + \alpha(v_1-v_2),
\end{cases}
\end{equation*}
where $u$ is the competitor, $v_1$ and $v_2$ are two phenotypes of the species $v=v_1+v_2$ with $d_1<d_2$, $\alpha>0$
is the mutation rate between phenotypes and the intrinsic growth rate $a(x)$ is heterogeneous. Strong competition means that $h>1$ and
$k>1$.

A similar system but with blind competition between $u$ and $v$ ($h=k=1$) is studied in a very recent work
by Cantrell, Cosner and Yu \cite{Cantrell_Cosner_Yu_2019}. In such a case, the long-time outcome still respects ``Unity is strength'', as can be expected. But what about the strong
competition case? The answer is likely quite complex and specific regimes should be considered instead of the general case;
numerical experiments could also be enlightening.
Note that even the existence of a pulsating traveling wave describing the confrontation
between $u$ and $v$ in a spatially periodic environment $a(x)$ is a difficult problem -- it 
might be solvable when $\alpha$ is so large that the interaction between $v_1$ and $v_2$
is essentially cooperative \cite{Cantrell_Cosner_Yu_2018}, but small values of $\alpha$ are elusive.

\subsection{Comparison between the Lotka--Volterra competition and other competition models}
Depending on the exact biological problem, there might be a more suitable competition model than the Lotka--Volterra one, namely quadratic
competition with a constant rate. 
As mentioned for instance by Perthame \cite[Chapter 4, Section 4.10, Exercise p. 83]{Perthame_2015}, 
a very straightforward calculation shows that, for the system with cubic competition
\begin{equation*}
    \begin{cases}
    \partial_{t}u = \Delta u + u\left( 1-u \right)-huv^2,\\
    \partial_{t}v = d\Delta v + rv\left( 1-v \right)-ku^2v,
\end{cases}
\end{equation*}
the speed of the traveling wave connecting $(1,0)$ and $(0,1)$ has the sign of $k-rh$. Surprisingly, the value of the diffusion rate $d$ does not
matter at all. This definitely shows that the conclusions drawn from the Lotka--Volterra case cannot be generalized without care.
Specific competition models require specific studies.

Note that the cubic competition model above is known in the physics literature as the Gross--Pitaevskii model (\textit{e.g.}, 
\cite{Dancer_Wang_Zh,Soave_Zilio_2015} and references therein) and in the mathematical biology literature as the 
Gilpin--Ayala model \cite{Gilpin_Ayala}. For other alternatives, refer for instance to Schoener
\cite{Schoener_1976} and Van Vuuren and Norbury \cite{Van_Vuuren_Norbury_2000} and references therein.

\section{Conclusion}
All existing results on strongly competing systems in homogeneous environments and
with populations differing only in diffusion rate concur: ``Disunity is strength''. In other words,
the fast dispersers win and chase the slow dispersers. Although the map of rigorous results 
(\figref{homogeneous_env}) is far from being complete, numerical investigations (\figref{numerics}) 
tend to confirm this is always the case: 
whatever the values of the two parameters (the interpopulation competition rate and the ratio 
between the dispersal rates) are, fast dispersers prevail. When the two populations differ 
also in intrinsic growth rate
or in interpopulation competition rate, the picture is more complicated but numerically
the conclusion seems to remain true close enough to symmetry (\figref{ADD_2012}).

In heterogeneous environments, the picture is less clear, in particular due to the possibilities of blocking or repelling
phenomena. Existing results concern the case of very strongly competing populations
in one-dimensional spatially periodic, temporally homogeneous environments.
On one hand, numerical investigations indicate that uniform dispersal rates are preferred over
highly oscillating ones. On the other hand, analytical works on uniform dispersal rates with periodic growth and competition
rates prove that a very large dispersal rate remains a definitive competitive advantage.
This last conclusion is surprisingly different from the ``Unity is strength''-type one obtained in 1998 \cite{Dockery_1998} 
for spatially heterogeneous environments where the competitive pressure is exerted blindly by individuals instead of
being mainly targeted at individuals of the other population.

More general heterogeneities, different dispersal strategies and different competition models
form very interesting but completely open problems that need more
attention from the community.

\section*{Acknowledgments}
The author thanks Florence D\'{e}barre for the attention she paid to this work and two anonymous referees for very valuable
suggestions.

\appendix
\section{On the numerical simulation for the wave speed in homogeneous environments}

First, for numerical convenience, change the spatial variable 
$x$ into $x\sqrt{d}$, so that the diffusion rate of $u$ becomes $1/d<1$, that of $v$ becomes $1$ and the wave speed $c_{k,d}$ 
becomes $c_{k,d}/\sqrt{d}$. 

Then, for any couple $(k,d)\in[1,21]^2$ (with steps of size $.1$), run a standard semi-implicit finite difference scheme 
in a bounded 1D domain of size $40\gg\max\left( 1,1/d \right)$ (steps of size $.02$) and during a time equal to $40$ (steps of size 
$.02$). 

Departing from piecewise-constant wave-like initial data, the numerical solution rapidly converges to the traveling wave. 
The speed $c_{k,d}/\sqrt{d}$ is evaluated by tracking the motion of the $.5$-level set of $v$ between $t=32$ and $t=40$.

\bibliographystyle{plain}
\bibliography{ref}

\end{document}

%% file: bistable_competitive_front.pdf_tex
\begingroup%
  \makeatletter%
  \providecommand\color[2][]{%
    \errmessage{(Inkscape) Color is used for the text in Inkscape, but the package 'color.sty' is not loaded}%
    \renewcommand\color[2][]{}%
  }%
  \providecommand\transparent[1]{%
    \errmessage{(Inkscape) Transparency is used (non-zero) for the text in Inkscape, but the package 'transparent.sty' is not loaded}%
    \renewcommand\transparent[1]{}%
  }%
  \providecommand\rotatebox[2]{#2}%
  \newcommand*\fsize{\dimexpr\f@size pt\relax}%
  \newcommand*\lineheight[1]{\fontsize{\fsize}{#1\fsize}\selectfont}%
  \ifx\svgwidth\undefined%
    \setlength{\unitlength}{302.78189111bp}%
    \ifx\svgscale\undefined%
      \relax%
    \else%
      \setlength{\unitlength}{\unitlength * \real{\svgscale}}%
    \fi%
  \else%
    \setlength{\unitlength}{\svgwidth}%
  \fi%
  \global\let\svgwidth\undefined%
  \global\let\svgscale\undefined%
  \makeatother%
  \begin{picture}(1,0.47865673)%
    \lineheight{1}%
    \setlength\tabcolsep{0pt}%
    \put(0,0){\includegraphics[width=\unitlength,page=1]{bistable_competitive_front.pdf}}%
    \put(0.29831481,0.30210723){\color[rgb]{0.49019608,0.49019608,0.49019608}\makebox(0,0)[lt]{\begin{minipage}{0.1462488\unitlength}\raggedright $u(t,x)$\end{minipage}}}%
    \put(0.68322468,0.08487188){\color[rgb]{0,0,0}\makebox(0,0)[lt]{\begin{minipage}{0.17111956\unitlength}\raggedright $u(t\pm1,x)$\end{minipage}}}%
    \put(0.52542661,0.19083139){\color[rgb]{1,0,0}\makebox(0,0)[lt]{\begin{minipage}{0.04706951\unitlength}\raggedright $c$\end{minipage}}}%
    \put(0.97635314,0.00856724){\color[rgb]{0,0,0}\makebox(0,0)[lt]{\lineheight{1.25}\smash{\begin{tabular}[t]{l}$x$\end{tabular}}}}%
    \put(0.03147781,0.46857327){\color[rgb]{0,0,0}\makebox(0,0)[lt]{\lineheight{1.25}\smash{\begin{tabular}[t]{l}$u,v$\end{tabular}}}}%
    \put(-0.00110166,0.34625125){\color[rgb]{0,0,0}\makebox(0,0)[lt]{\lineheight{1.25}\smash{\begin{tabular}[t]{l}$1$\end{tabular}}}}%
    \put(-0.00110166,0.00930744){\color[rgb]{0,0,0}\makebox(0,0)[lt]{\lineheight{1.25}\smash{\begin{tabular}[t]{l}$0$\end{tabular}}}}%
    \put(0,0){\includegraphics[width=\unitlength,page=2]{bistable_competitive_front.pdf}}%
    \put(0.31130357,0.07594996){\color[rgb]{0.49019608,0.49019608,0.49019608}\makebox(0,0)[lt]{\begin{minipage}{0.1462488\unitlength}\raggedright $v(t,x)$\end{minipage}}}%
    \put(0.73747186,0.30415275){\color[rgb]{0,0,0}\makebox(0,0)[lt]{\begin{minipage}{0.17111956\unitlength}\raggedright $v(t\pm1,x)$\end{minipage}}}%
  \end{picture}%
\endgroup%

%% file: UINS_with_references.pdf_tex
\begingroup%
  \makeatletter%
  \providecommand\color[2][]{%
    \errmessage{(Inkscape) Color is used for the text in Inkscape, but the package 'color.sty' is not loaded}%
    \renewcommand\color[2][]{}%
  }%
  \providecommand\transparent[1]{%
    \errmessage{(Inkscape) Transparency is used (non-zero) for the text in Inkscape, but the package 'transparent.sty' is not loaded}%
    \renewcommand\transparent[1]{}%
  }%
  \providecommand\rotatebox[2]{#2}%
  \newcommand*\fsize{\dimexpr\f@size pt\relax}%
  \newcommand*\lineheight[1]{\fontsize{\fsize}{#1\fsize}\selectfont}%
  \ifx\svgwidth\undefined%
    \setlength{\unitlength}{415.99631963bp}%
    \ifx\svgscale\undefined%
      \relax%
    \else%
      \setlength{\unitlength}{\unitlength * \real{\svgscale}}%
    \fi%
  \else%
    \setlength{\unitlength}{\svgwidth}%
  \fi%
  \global\let\svgwidth\undefined%
  \global\let\svgscale\undefined%
  \makeatother%
  \begin{picture}(1,0.89578128)%
    \lineheight{1}%
    \setlength\tabcolsep{0pt}%
    \put(0,0){\includegraphics[width=\unitlength,page=1]{UINS_with_references.pdf}}%
    \put(0.03494022,0.0087927){\color[rgb]{0,0,0}\makebox(0,0)[lt]{\lineheight{1.25}\smash{\begin{tabular}[t]{l}$(1,1)$\end{tabular}}}}%
    \put(0.95911046,0.00250949){\color[rgb]{0,0,0}\makebox(0,0)[lt]{\lineheight{1.25}\smash{\begin{tabular}[t]{l}$d$\end{tabular}}}}%
    \put(0.04768865,0.87295417){\color[rgb]{0,0,0}\makebox(0,0)[lt]{\lineheight{1.25}\smash{\begin{tabular}[t]{l}$k$\end{tabular}}}}%
    \put(0.38092331,0.00695307){\color[rgb]{0,0,0}\makebox(0,0)[lt]{\lineheight{1.25}\smash{\begin{tabular}[t]{l}$4$\end{tabular}}}}%
    \put(0.50572161,0.00423963){\color[rgb]{0,0,0}\makebox(0,0)[lt]{\lineheight{1.25}\smash{\begin{tabular}[t]{l}$5.5$\end{tabular}}}}%
    \put(0.034867,0.11340205){\color[rgb]{0,0,0}\makebox(0,0)[lt]{\lineheight{1.25}\smash{\begin{tabular}[t]{l}$1.25$\end{tabular}}}}%
    \put(0.00211278,0.14809641){\color[rgb]{0,0,0}\makebox(0,0)[lt]{\lineheight{1.25}\smash{\begin{tabular}[t]{l}$1.33\dots$\end{tabular}}}}%
    \put(0.00397975,0.35846389){\color[rgb]{0,0,0}\makebox(0,0)[lt]{\lineheight{1.25}\smash{\begin{tabular}[t]{l}$1.83\dots$\end{tabular}}}}%
    \put(0,0){\includegraphics[width=\unitlength,page=2]{UINS_with_references.pdf}}%
    \put(0.90456481,0.51318971){\color[rgb]{0,0,0}\rotatebox{-90}{\makebox(0,0)[lt]{\lineheight{1.25}\smash{\begin{tabular}[t]{l}\cite{Alzahrani_Davidson_Dodds_2010} ($d>\underline{d}(k)$)\end{tabular}}}}}%
    \put(0.47105455,0.79932668){\color[rgb]{0,0,0}\makebox(0,0)[lt]{\lineheight{1.25}\smash{\begin{tabular}[t]{l}\cite{Girardin_Nadin_2015} ($k>\underline{k}(d)$)\end{tabular}}}}%
    \put(0.06129862,0.43540052){\color[rgb]{0,0,0}\makebox(0,0)[lt]{\lineheight{1.25}\smash{\begin{tabular}[t]{l}$2$\end{tabular}}}}%
    \put(0.00270435,0.27584376){\color[rgb]{0,0,0}\makebox(0,0)[lt]{\lineheight{1.25}\smash{\begin{tabular}[t]{l}$1.66\dots$\end{tabular}}}}%
    \put(0.57185534,0.0057821){\color[rgb]{0,0,0}\makebox(0,0)[lt]{\lineheight{1.25}\smash{\begin{tabular}[t]{l}$6$\end{tabular}}}}%
    \put(0.4678675,0.00509951){\color[rgb]{0,0,0}\makebox(0,0)[lt]{\lineheight{1.25}\smash{\begin{tabular}[t]{l}$5$\end{tabular}}}}%
    \put(0.40499641,0.42176154){\color[rgb]{0,0,0}\makebox(0,0)[lt]{\lineheight{1.25}\smash{\begin{tabular}[t]{l}\cite{Ma_Huang_Ou_2019} ($d(k-1)=4$)\end{tabular}}}}%
    \put(0.5316096,0.35313157){\color[rgb]{0,0,0}\makebox(0,0)[lt]{\lineheight{1.25}\smash{\begin{tabular}[t]{l}\cite{Rodrigo_Mimura_2001}\end{tabular}}}}%
    \put(0.40157563,0.12875918){\color[rgb]{0,0,0}\makebox(0,0)[lt]{\lineheight{1.25}\smash{\begin{tabular}[t]{l}\cite{Guo_Lin_2013}\end{tabular}}}}%
    \put(0.10198746,0.04716918){\color[rgb]{0,0,0}\makebox(0,0)[lt]{\lineheight{1.25}\smash{\begin{tabular}[t]{l}\cite{Risler_2017}\end{tabular}}}}%
  \end{picture}%
\endgroup%

%% file: ADD_remake.pdf_tex
\begingroup%
  \makeatletter%
  \providecommand\color[2][]{%
    \errmessage{(Inkscape) Color is used for the text in Inkscape, but the package 'color.sty' is not loaded}%
    \renewcommand\color[2][]{}%
  }%
  \providecommand\transparent[1]{%
    \errmessage{(Inkscape) Transparency is used (non-zero) for the text in Inkscape, but the package 'transparent.sty' is not loaded}%
    \renewcommand\transparent[1]{}%
  }%
  \providecommand\rotatebox[2]{#2}%
  \newcommand*\fsize{\dimexpr\f@size pt\relax}%
  \newcommand*\lineheight[1]{\fontsize{\fsize}{#1\fsize}\selectfont}%
  \ifx\svgwidth\undefined%
    \setlength{\unitlength}{433.83059044bp}%
    \ifx\svgscale\undefined%
      \relax%
    \else%
      \setlength{\unitlength}{\unitlength * \real{\svgscale}}%
    \fi%
  \else%
    \setlength{\unitlength}{\svgwidth}%
  \fi%
  \global\let\svgwidth\undefined%
  \global\let\svgscale\undefined%
  \makeatother%
  \begin{picture}(1,0.62668458)%
    \lineheight{1}%
    \setlength\tabcolsep{0pt}%
    \put(0,0){\includegraphics[width=\unitlength,page=1]{ADD_remake.pdf}}%
    \put(0.02979555,0.02029698){\color[rgb]{0,0,0}\makebox(0,0)[lt]{\lineheight{1.25}\smash{\begin{tabular}[t]{l}$(0,0)$\end{tabular}}}}%
    \put(0.91597423,0.01427206){\color[rgb]{0,0,0}\makebox(0,0)[lt]{\lineheight{1.25}\smash{\begin{tabular}[t]{l}$h$\end{tabular}}}}%
    \put(0.04568721,0.5946671){\color[rgb]{0,0,0}\makebox(0,0)[lt]{\lineheight{1.25}\smash{\begin{tabular}[t]{l}$d$\end{tabular}}}}%
    \put(0.26889487,0.00697606){\color[rgb]{0,0,0}\makebox(0,0)[lt]{\lineheight{1.25}\smash{\begin{tabular}[t]{l}$\sqrt{\frac{k}{r}}$\end{tabular}}}}%
    \put(0.43759106,0.00575362){\color[rgb]{0,0,0}\makebox(0,0)[lt]{\lineheight{1.25}\smash{\begin{tabular}[t]{l}$\frac{k}{r}$\end{tabular}}}}%
    \put(0.64173783,0.00575362){\color[rgb]{0,0,0}\makebox(0,0)[lt]{\lineheight{1.25}\smash{\begin{tabular}[t]{l}$\left(\frac{k}{r}\right)^2$\end{tabular}}}}%
    \put(0.03606429,0.17825656){\color[rgb]{0,0,0}\makebox(0,0)[lt]{\lineheight{1.25}\smash{\begin{tabular}[t]{l}$1$\end{tabular}}}}%
    \put(0,0){\includegraphics[width=\unitlength,page=2]{ADD_remake.pdf}}%
    \put(0.20655066,0.22579211){\color[rgb]{0,0,0}\makebox(0,0)[lt]{\lineheight{1.25}\smash{\begin{tabular}[t]{l}$c_{d,r,h,k}>0$\end{tabular}}}}%
    \put(0.58428335,0.49595034){\color[rgb]{0,0,0}\makebox(0,0)[lt]{\lineheight{1.25}\smash{\begin{tabular}[t]{l}$c_{d,r,h,k}<0$\end{tabular}}}}%
    \put(0.33001668,0.40671254){\color[rgb]{0,0,0}\makebox(0,0)[lt]{\lineheight{1.25}\smash{\begin{tabular}[t]{l}$c_{d,r,h,k}=0$\end{tabular}}}}%
  \end{picture}%
\endgroup%

%% file: periodic_coexistence_steady_state.pdf_tex
\begingroup%
  \makeatletter%
  \providecommand\color[2][]{%
    \errmessage{(Inkscape) Color is used for the text in Inkscape, but the package 'color.sty' is not loaded}%
    \renewcommand\color[2][]{}%
  }%
  \providecommand\transparent[1]{%
    \errmessage{(Inkscape) Transparency is used (non-zero) for the text in Inkscape, but the package 'transparent.sty' is not loaded}%
    \renewcommand\transparent[1]{}%
  }%
  \providecommand\rotatebox[2]{#2}%
  \newcommand*\fsize{\dimexpr\f@size pt\relax}%
  \newcommand*\lineheight[1]{\fontsize{\fsize}{#1\fsize}\selectfont}%
  \ifx\svgwidth\undefined%
    \setlength{\unitlength}{302.78189111bp}%
    \ifx\svgscale\undefined%
      \relax%
    \else%
      \setlength{\unitlength}{\unitlength * \real{\svgscale}}%
    \fi%
  \else%
    \setlength{\unitlength}{\svgwidth}%
  \fi%
  \global\let\svgwidth\undefined%
  \global\let\svgscale\undefined%
  \makeatother%
  \begin{picture}(1,0.47865673)%
    \lineheight{1}%
    \setlength\tabcolsep{0pt}%
    \put(0,0){\includegraphics[width=\unitlength,page=1]{periodic_coexistence_steady_state.pdf}}%
    \put(0.97635314,0.00856724){\color[rgb]{0,0,0}\makebox(0,0)[lt]{\lineheight{1.25}\smash{\begin{tabular}[t]{l}$x$\end{tabular}}}}%
    \put(0.03147781,0.46857327){\color[rgb]{0,0,0}\makebox(0,0)[lt]{\lineheight{1.25}\smash{\begin{tabular}[t]{l}$u,v$\end{tabular}}}}%
    \put(-0.00110166,0.34319507){\color[rgb]{1,0,0}\makebox(0,0)[lt]{\lineheight{1.25}\smash{\begin{tabular}[t]{l}$1$\end{tabular}}}}%
    \put(-0.00110166,0.00930744){\color[rgb]{0,0,0}\makebox(0,0)[lt]{\lineheight{1.25}\smash{\begin{tabular}[t]{l}$0$\end{tabular}}}}%
    \put(0.03246998,0.44707417){\color[rgb]{1,0,0}\makebox(0,0)[lt]{\lineheight{1.25}\smash{\begin{tabular}[t]{l}$\mu$\end{tabular}}}}%
    \put(0,0){\includegraphics[width=\unitlength,page=2]{periodic_coexistence_steady_state.pdf}}%
    \put(0.02483146,-0.02368539){\color[rgb]{0,0,0}\makebox(0,0)[lt]{\lineheight{1.25}\smash{\begin{tabular}[t]{l}$0$\end{tabular}}}}%
    \put(0.91799981,-0.02483137){\color[rgb]{0,0,0}\makebox(0,0)[lt]{\lineheight{1.25}\smash{\begin{tabular}[t]{l}$L$\end{tabular}}}}%
    \put(0.23150559,0.14287638){\color[rgb]{0,0,0}\makebox(0,0)[lt]{\lineheight{1.25}\smash{\begin{tabular}[t]{l}$v$\end{tabular}}}}%
    \put(0.68305616,0.14440449){\color[rgb]{0,0,0}\makebox(0,0)[lt]{\lineheight{1.25}\smash{\begin{tabular}[t]{l}$u$\end{tabular}}}}%
  \end{picture}%
\endgroup%

%% file: heatmap.tex
\setlength{\unitlength}{1pt}
\begin{picture}(0,0)
\includegraphics{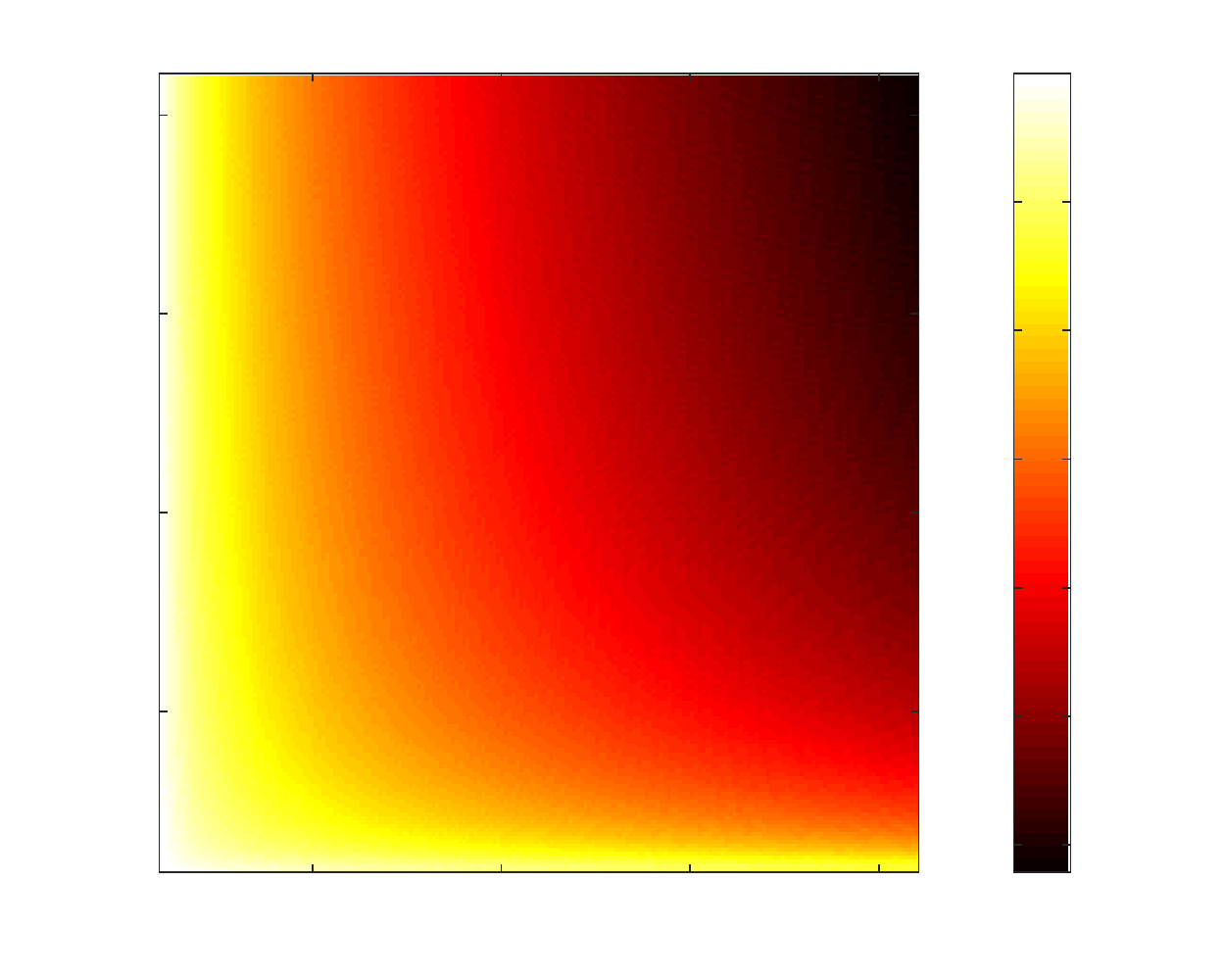}
\end{picture}%
\begin{picture}(360,288)(0,0)
\fontsize{10}{0}
\selectfont\put(158.4,15.686){\makebox(0,0)[t]{\textcolor[rgb]{0.15,0.15,0.15}{{$d$}}}}
\fontsize{10}{0}
\selectfont\put(25.7955,149.04){\rotatebox{90}{\makebox(0,0)[b]{\textcolor[rgb]{0.15,0.15,0.15}{{$k$}}}}}
\fontsize{10}{0}
\selectfont\put(91.7748,26.686){\makebox(0,0)[t]{\textcolor[rgb]{0.15,0.15,0.15}{{5}}}}
\fontsize{10}{0}
\selectfont\put(147.296,26.686){\makebox(0,0)[t]{\textcolor[rgb]{0.15,0.15,0.15}{{10}}}}
\fontsize{10}{0}
\selectfont\put(202.817,26.686){\makebox(0,0)[t]{\textcolor[rgb]{0.15,0.15,0.15}{{15}}}}
\fontsize{10}{0}
\selectfont\put(258.338,26.686){\makebox(0,0)[t]{\textcolor[rgb]{0.15,0.15,0.15}{{20}}}}
\fontsize{10}{0}
\selectfont\put(41.7955,78.9761){\makebox(0,0)[r]{\textcolor[rgb]{0.15,0.15,0.15}{{5}}}}
\fontsize{10}{0}
\selectfont\put(41.7955,137.363){\makebox(0,0)[r]{\textcolor[rgb]{0.15,0.15,0.15}{{10}}}}
\fontsize{10}{0}
\selectfont\put(41.7955,195.749){\makebox(0,0)[r]{\textcolor[rgb]{0.15,0.15,0.15}{{15}}}}
\fontsize{10}{0}
\selectfont\put(41.7955,254.136){\makebox(0,0)[r]{\textcolor[rgb]{0.15,0.15,0.15}{{20}}}}
\fontsize{10}{0}
\selectfont\put(319.564,39.7394){\makebox(0,0)[l]{\textcolor[rgb]{0.15,0.15,0.15}{{-1.2}}}}
\fontsize{10}{0}
\selectfont\put(319.564,77.5162){\makebox(0,0)[l]{\textcolor[rgb]{0.15,0.15,0.15}{{-1}}}}
\fontsize{10}{0}
\selectfont\put(319.564,115.293){\makebox(0,0)[l]{\textcolor[rgb]{0.15,0.15,0.15}{{-0.8}}}}
\fontsize{10}{0}
\selectfont\put(319.564,153.07){\makebox(0,0)[l]{\textcolor[rgb]{0.15,0.15,0.15}{{-0.6}}}}
\fontsize{10}{0}
\selectfont\put(319.564,190.846){\makebox(0,0)[l]{\textcolor[rgb]{0.15,0.15,0.15}{{-0.4}}}}
\fontsize{10}{0}
\selectfont\put(319.564,228.623){\makebox(0,0)[l]{\textcolor[rgb]{0.15,0.15,0.15}{{-0.2}}}}
\fontsize{10}{0}
\selectfont\put(319.564,266.4){\makebox(0,0)[l]{\textcolor[rgb]{0.15,0.15,0.15}{{0}}}}
\end{picture}

%% file: level_sets.tex
\setlength{\unitlength}{1pt}
\begin{picture}(0,0)
\includegraphics{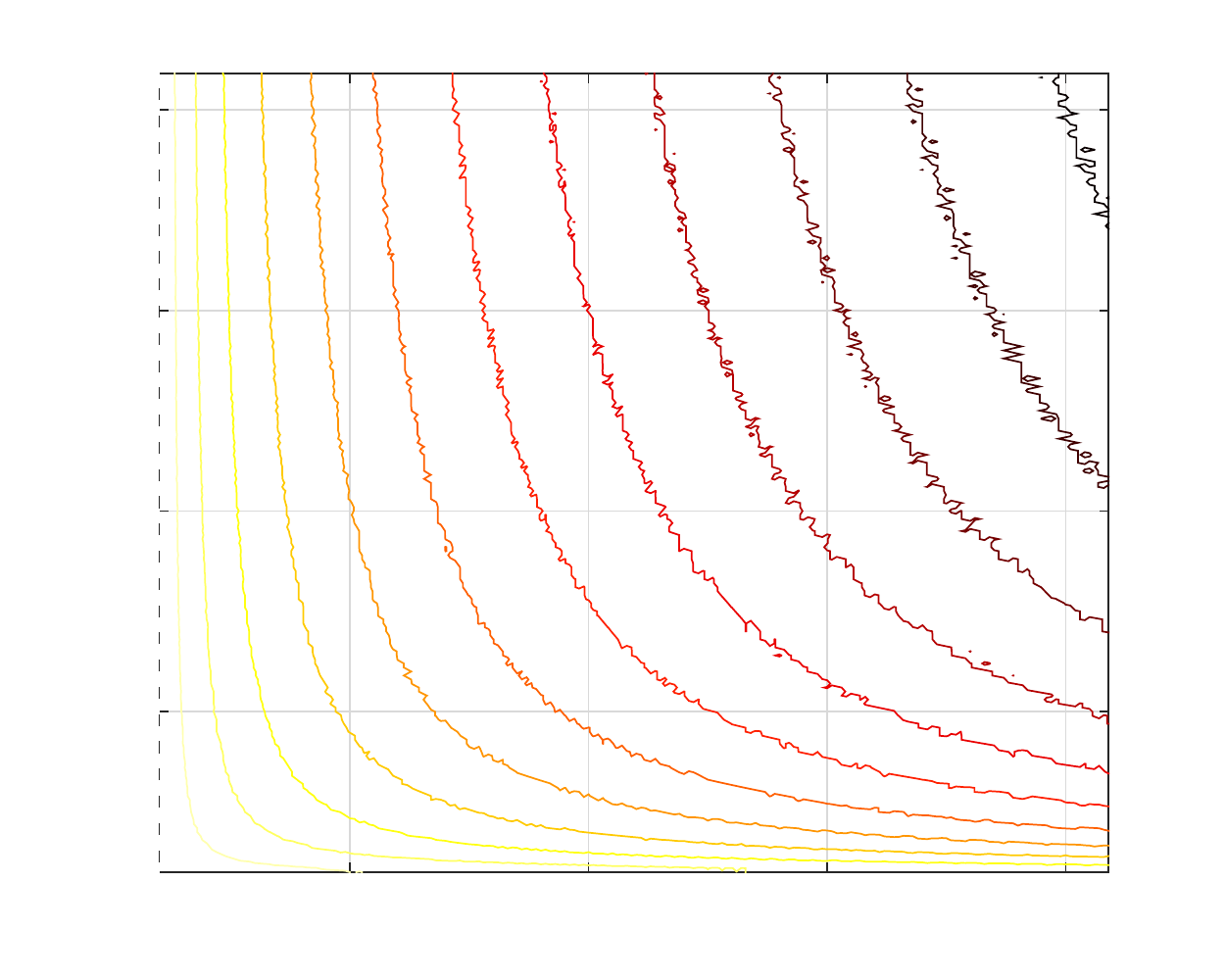}
\end{picture}%
\begin{picture}(360,288)(0,0)
\fontsize{10}{0}
\selectfont\put(102.88,26.686){\makebox(0,0)[t]{\textcolor[rgb]{0.15,0.15,0.15}{{5}}}}
\fontsize{10}{0}
\selectfont\put(172.981,26.686){\makebox(0,0)[t]{\textcolor[rgb]{0.15,0.15,0.15}{{10}}}}
\fontsize{10}{0}
\selectfont\put(243.081,26.686){\makebox(0,0)[t]{\textcolor[rgb]{0.15,0.15,0.15}{{15}}}}
\fontsize{10}{0}
\selectfont\put(313.182,26.686){\makebox(0,0)[t]{\textcolor[rgb]{0.15,0.15,0.15}{{20}}}}
\fontsize{10}{0}
\selectfont\put(41.8,78.8599){\makebox(0,0)[r]{\textcolor[rgb]{0.15,0.15,0.15}{{5}}}}
\fontsize{10}{0}
\selectfont\put(41.8,137.835){\makebox(0,0)[r]{\textcolor[rgb]{0.15,0.15,0.15}{{10}}}}
\fontsize{10}{0}
\selectfont\put(41.8,196.81){\makebox(0,0)[r]{\textcolor[rgb]{0.15,0.15,0.15}{{15}}}}
\fontsize{10}{0}
\selectfont\put(41.8,255.785){\makebox(0,0)[r]{\textcolor[rgb]{0.15,0.15,0.15}{{20}}}}
\fontsize{10}{0}
\selectfont\put(186.3,15.686){\makebox(0,0)[t]{\textcolor[rgb]{0.15,0.15,0.15}{{$d$}}}}
\fontsize{10}{0}
\selectfont\put(25.8,149.04){\rotatebox{90}{\makebox(0,0)[b]{\textcolor[rgb]{0.15,0.15,0.15}{{$k$}}}}}
\end{picture}

%% file: LV_bistable_survey.bbl
\begin{thebibliography}{10}

\bibitem{Alzahrani_Davidson_Dodds_2010}
Ebraheem~O. Alzahrani, Fordyce~A. Davidson, and Niall Dodds.
\newblock Travelling waves in near-degenerate bistable competition models.
\newblock {\em Math. Model. Nat. Phenom.}, 5(5):13--35, 2010.

\bibitem{Alzahrani_Davidson_Dodds_2012}
Ebraheem~O. Alzahrani, Fordyce~A. Davidson, and Niall Dodds.
\newblock Reversing invasion in bistable systems.
\newblock {\em J. Math. Biol.}, 65(6-7):1101--1124, 2012.

\bibitem{Amarasekare_20}
Priyanga Amarasekare.
\newblock Competitive coexistence in spatially structured environments: a
  synthesis.
\newblock {\em Ecology Letters}, 6(12):1109--1122, 2003.

\bibitem{Bao_Wang_2013}
Xiongxiong Bao and Zhi-Cheng Wang.
\newblock Existence and stability of time periodic traveling waves for a
  periodic bistable {L}otka--{V}olterra competition system.
\newblock {\em Journal of Differential Equations}, 255(8):2402--2435, 2013.

\bibitem{Barton_Turelli}
Nick~H. Barton and Michael Turelli.
\newblock Spatial waves of advance with bistable dynamics: Cytoplasmic and
  genetic analogues of {A}llee effects.
\newblock {\em The American Naturalist}, 178(3):E48--E75, 2011.
\newblock PMID: 21828986.

\bibitem{Benichou_2011}
Olivier B{\'{e}}nichou, Claude Loverdo, Michel Moreau, and Raphael Voituriez.
\newblock Intermittent search strategies.
\newblock {\em Rev. Mod. Phys.}, 83(1):81--129, Mar 2011.

\bibitem{Berestycki_Bouhours_Chapuisat}
Henri Berestycki, Juliette Bouhours, and Guillemette Chapuisat.
\newblock Front blocking and propagation in cylinders with varying cross
  section.
\newblock {\em Calc. Var. Partial Differential Equations}, 55(3):--44, 2016.

\bibitem{Berestycki_Ham_3}
Henri Berestycki and Fran{\c{c}}ois Hamel.
\newblock Front propagation in periodic excitable media.
\newblock {\em Comm. Pure Appl. Math.}, 55(8):949--1032, 2002.

\bibitem{Berestycki_Ham-2}
Henri Berestycki and Fran{\c{c}}ois Hamel.
\newblock Generalized travelling waves for reaction-diffusion equations.
\newblock In {\em Perspectives in nonlinear partial differential equations},
  volume 446, chapter Contemp. Math., pages 101--123. Amer. Math. Soc.,
  Providence, RI, 2007.

\bibitem{Berestycki_Ham-1}
Henri Berestycki and Fran{\c{c}}ois Hamel.
\newblock Generalized transition waves and their properties.
\newblock {\em Comm. Pure Appl. Math.}, 65(5):592--648, 2012.

\bibitem{Berestycki_Ham_1}
Henri Berestycki, Fran{\c{c}}ois Hamel, and Lionel Roques.
\newblock Analysis of the periodically fragmented environment model. {I}.
  {S}pecies persistence.
\newblock {\em J. Math. Biol.}, 51(1):75--113, 2005.

\bibitem{Berestycki_Ham_2}
Henri Berestycki, Fran{\c{c}}ois Hamel, and Lionel Roques.
\newblock Analysis of the periodically fragmented environment model. {II}.
  {B}iological invasions and pulsating travelling fronts.
\newblock {\em J. Math. Pures Appl. (9)}, 84(8):1101--1146, 2005.

\bibitem{Berestycki_Zilio_17}
Henri Berestycki and Alessandro Zilio.
\newblock Predator-prey models with competition: The emergence of
  territoriality.
\newblock {\em The American Naturalist}, 193(3):436--446, 2019.
\newblock PMID: 30794454.

\bibitem{Boeye_2014}
Jeroen Boeye, Alexander Kubisch, and Dries Bonte.
\newblock Habitat structure mediates spatial segregation and therefore
  coexistence.
\newblock {\em Landscape Ecology}, 29(4):593--604, Apr 2014.

\bibitem{Bohn_Amundsen_}
Thomas B{\o}hn, Per-Arne Amundsen, and Ashley Sparrow.
\newblock Competitive exclusion after invasion?
\newblock {\em Biological Invasions}, 10(3):359--368, Mar 2008.

\bibitem{Brown_71}
James~H. Brown.
\newblock Mechanisms of competitive exclusion between two species of chipmunks.
\newblock {\em Ecology}, 52(2):305--311, 2011.

\bibitem{Cantrell_Cosner_Yu_2018}
Robert~Stephen Cantrell, Chris Cosner, and Xiao Yu.
\newblock Dynamics of populations with individual variation in dispersal on
  bounded domains.
\newblock {\em Journal of biological dynamics}, 12(1):288--317, 2018.

\bibitem{Cantrell_Cosner_Yu_2019}
Robert~Stephen Cantrell, Chris Cosner, and Xiao Yu.
\newblock Populations with individual variation in dispersal in heterogeneous
  environments: Dynamics and competition with simply diffusing populations.
\newblock {\em Science China Mathematics}, Jan 2020.

\bibitem{Cantrell_Cosner_Ruan_2009}
Stephen Cantrell, Chris Cosner, and Shigui Ruan.
\newblock {\em Spatial ecology}.
\newblock CRC Press, 2009.

\bibitem{Carrere_2017_2}
C\'{e}cile Carr\`ere.
\newblock Optimization of an \textit{in vitro} chemotherapy to avoid resistant
  tumours.
\newblock {\em J. Theoret. Biol.}, 413:24--33, 2017.

\bibitem{Carrere_2017}
C{\'{e}}cile Carr{\`{e}}re.
\newblock Spreading speeds for a two-species competition-diffusion system.
\newblock {\em J. Differential Equations}, 264(3):2133--2156, 2018.

\bibitem{Dancer_Wang_Zh}
Edward~N. Dancer, Kelei Wang, and Zhitao Zhang.
\newblock The limit equation for the {G}ross--{P}itaevskii equations and {S}.
  {T}erracini's conjecture.
\newblock {\em Journal of Functional Analysis}, 262(3):1087--1131, 2012.

\bibitem{Deforet_2017}
Maxime Deforet, Carlos Carmona-Fontaine, Kirill~S. Korolev, and Joao~B. Xavier.
\newblock Evolution at the edge of expanding populations.
\newblock {\em The American Naturalist}, 194(3):291--305, 2019.

\bibitem{Dockery_1998}
Jack Dockery, Vivian Hutson, Konstantin Mischaikow, and Mark Pernarowski.
\newblock The evolution of slow dispersal rates: a reaction diffusion model.
\newblock {\em J. Math. Biol.}, 37(1):61--83, 1998.

\bibitem{Dreher_2007}
Michael Dreher.
\newblock Analysis of a population model with strong cross-diffusion in
  unbounded domains.
\newblock {\em Proc. Roy. Soc. Edinburgh Sect. A}, 138(4):769--786, 2008.

\bibitem{Du_Li_Wu_2018}
Li-Jun Du, Wan-Tong Li, and Shi-Liang Wu.
\newblock Propagation phenomena for a bistable {L}otka-{V}olterra competition
  system with advection in a periodic habitat.
\newblock {\em Zeitschrift f{\"u}r angewandte Mathematik und Physik}, 71(1):11,
  2019.

\bibitem{Debarre_Lenormand_11}
Florence Débarre and Thomas Lenormand.
\newblock Distance-limited dispersal promotes coexistence at habitat
  boundaries: reconsidering the competitive exclusion principle.
\newblock {\em Ecology Letters}, 14(3):260--266, 2011.

\bibitem{Octave}
John~W. Eaton, David Bateman, S{\o}ren Hauberg, and Rik Wehbring.
\newblock {\em {GNU Octave} version 5.1.0 manual: a high-level interactive
  language for numerical computations}, 2019.

\bibitem{Fisher_1937}
Ronald~Aylmer Fisher.
\newblock The wave of advance of advantageous genes.
\newblock {\em Annals of eugenics}, 7(4):355--369, 1937.

\bibitem{Gardner_1982}
Robert~A. Gardner.
\newblock Existence and stability of travelling wave solutions of competition
  models: a degree theoretic approach.
\newblock {\em J. Differential Equations}, 44(3):343--364, 1982.

\bibitem{Gartner_Freidlin}
J\"{u}rgen G\"{a}rtner and Mark~I. Freidlin.
\newblock The propagation of concentration waves in periodic and random media.
\newblock {\em Dokl. Akad. Nauk SSSR}, 249(3):521--525, 1979.

\bibitem{Gause_1934}
Georgii~F. Gause.
\newblock {\em The {S}truggle for {E}xistence}.
\newblock Baltimore, 1934.

\bibitem{Gilpin_Ayala}
Michael~E. Gilpin and Francisco~J. Ayala.
\newblock Global models of growth and competition.
\newblock {\em Proceedings of the National Academy of Sciences},
  70(12):3590--3593, 1973.

\bibitem{Girardin_2016}
L\'{e}o Girardin.
\newblock Competition in periodic media: {I} -- {E}xistence of pulsating
  fronts.
\newblock {\em Discrete and Continuous Dynamical Systems - Series B},
  22(4):1341--1360, 2017.

\bibitem{Girardin_Nadin_2015}
L\'{e}o Girardin and Gr\'{e}goire Nadin.
\newblock Travelling waves for diffusive and strongly competitive systems:
  relative motility and invasion speed.
\newblock {\em European J. Appl. Math.}, 26(4):521--534, 2015.

\bibitem{Girardin_Nadin_2016}
L\'{e}o Girardin and Gr\'{e}goire Nadin.
\newblock Competition in periodic media: {II} -- {S}egregative limit of
  pulsating fronts and {{``}Unity is not Strength{''}}-type result.
\newblock {\em Journal of Differential Equations}, 265(1):98--156, 2018.

\bibitem{Girardin_Zilio}
L\'{e}o Girardin and Alessandro Zilio.
\newblock Competition in periodic media: {III} -- {E}xistence \& stability of
  segregated periodic coexistence states.
\newblock {\em Journal of Dynamics and Differential Equations}, 2019.

\bibitem{Giuggioli_2010}
Luca Giuggioli and Frederic Bartumeus.
\newblock Animal movement, search strategies and behavioural ecology: a
  cross-disciplinary way forward.
\newblock {\em Journal of Animal Ecology}, 79:906{--}909, 2010.

\bibitem{Guo_Lin_2013}
Jong-Shenq Guo and Ying-Chih Lin.
\newblock The sign of the wave speed for the {L}otka-{V}olterra
  competition-diffusion system.
\newblock {\em Commun. Pure Appl. Anal.}, 12(5):2083--2090, 2013.

\bibitem{Heinze_Schweiz}
Steffen Heinze and Ben Schweizer.
\newblock Creeping fronts in degenerate reaction-diffusion systems.
\newblock {\em Nonlinearity}, 18(6):2455--2476, 2005.

\bibitem{Heinze_et_al_2004}
Steffen Heinze, Ben Schweizer, and Hartmut Schwetlick.
\newblock Existence of front solutions in degenerate reaction diffusion
  systems.
\newblock {\em Preprint}, 2004.

\bibitem{Hutridurga_Ven}
Harsha Hutridurga and Chandrasekhar Venkataraman.
\newblock Heterogeneity and strong competition in ecology.
\newblock {\em European Journal of Applied Mathematics}, page 1–25, 2018.

\bibitem{Hutson_Martinez_2003}
Vivian Hutson, Salom\'{e} Martinez, Konstantin Mischaikow, and Glenn~T.
  Vickers.
\newblock The evolution of dispersal.
\newblock {\em J. Math. Biol.}, 47(6):483--517, 2003.

\bibitem{Hutson_Mischaikow_Polacik}
Vivian Hutson, Konstantin Mischaikow, and Peter Pol{\'{a}}{\v{c}}ik.
\newblock The evolution of dispersal rates in a heterogeneous time-periodic
  environment.
\newblock {\em J. Math. Biol.}, 43(6):501--533, 2001.

\bibitem{Iannelli_Pugliese}
Mimmo Iannelli and Andrea Pugliese.
\newblock {\em An introduction to mathematical population dynamics}, volume~79
  of {\em Unitext}.
\newblock Springer, Cham, 2014.
\newblock Along the trail of Volterra and Lotka, La Matematica per il 3+2.

\bibitem{Kan_on_1995}
Yukio Kan-on.
\newblock Parameter dependence of propagation speed of travelling waves for
  competition-diffusion equations.
\newblock {\em SIAM J. Math. Anal.}, 26(2):340--363, 1995.

\bibitem{Kan-on_Yanagida_1993}
Yukio Kan-on and Eiji Yanagida.
\newblock Existence of nonconstant stable equilibria in competition-diffusion
  equations.
\newblock {\em Hiroshima Math. J.}, 23(1):193--221, 1993.

\bibitem{Kimura_1965}
Motoo Kimura.
\newblock A stochastic model concerning the maintenance of genetic variability
  in quantitative characters.
\newblock {\em Proceedings of the National Academy of Sciences of the United
  States of America}, 54(3):731--736, 1965.

\bibitem{Kishimoto_Wein}
Kazuo Kishimoto and Hans~F. Weinberger.
\newblock The spatial homogeneity of stable equilibria of some
  reaction-diffusion systems on convex domains.
\newblock {\em J. Differential Equations}, 58(1):15--21, 1985.

\bibitem{KPP_1937}
Andrei~N. Kolmogorov, I.~G. Petrovsky, and N.~S. Piskunov.
\newblock {\'Etude} de l'\'equation de la diffusion avec croissance de la
  quantit\'e de mati\`ere et son application \`a un probl\`eme biologique.
\newblock {\em Bulletin Universit\'e d'\'Etat {\`{a}} Moscou}, 1:1--25, 1937.

\bibitem{Lotka_1924}
Alfred~J. Lotka.
\newblock {\em Elements of {P}hysical {B}iology,}.
\newblock Baltimore, 1924.

\bibitem{Ma_Huang_Ou_2019}
Manjun Ma, Zhe Huang, and Chunhua Ou.
\newblock Speed of the traveling wave for the bistable {L}otka--{V}olterra
  competition model.
\newblock {\em Nonlinearity}, 32(9):3143--3162, jul 2019.

\bibitem{Maciel_Lutscher_18}
Gabriel~Andreguetto Maciel and Frithjof Lutscher.
\newblock Movement behaviour determines competitive outcome and spread rates in
  strongly heterogeneous landscapes.
\newblock {\em Theoretical Ecology}, 11(3):351--365, Sep 2018.

\bibitem{Matano_79}
Hiroshi Matano.
\newblock Asymptotic behavior and stability of solutions of semilinear
  diffusion equations.
\newblock {\em Publ. Res. Inst. Math. Sci.}, 15(2):401--454, 1979.

\bibitem{Matano_Mimura_83}
Hiroshi Matano and Masayasu Mimura.
\newblock Pattern formation in competition-diffusion systems in nonconvex
  domains.
\newblock {\em Publ. Res. Inst. Math. Sci.}, 19(3):1049--1079, 1983.

\bibitem{Melbourne_Corn_07}
Brett~A. Melbourne, Howard~V. Cornell, Kendi~F. Davies, Christopher~J. Dugaw,
  Sarah Elmendorf, Amy~L. Freestone, Richard~J. Hall, Susan Harrison, Alan
  Hastings, Matt Holland, Marcel Holyoak, John Lambrinos, Kara Moore, and
  Hiroyuki Yokomizo.
\newblock Invasion in a heterogeneous world: resistance, coexistence or hostile
  takeover?
\newblock {\em Ecology Letters}, 10(1):77--94, 2007.

\bibitem{Mitani_Watts_Amsler}
John~C. Mitani, David~P. Watts, and Sylvia~J. Amsler.
\newblock Lethal intergroup aggression leads to territorial expansion in wild
  chimpanzees.
\newblock {\em Current Biology}, 20(12):--507, 2010.

\bibitem{Ninomiya_1995}
Hirokazu Ninomiya.
\newblock Separatrices of competition-diffusion equations.
\newblock {\em J. Math. Kyoto Univ.}, 35(3):539--567, 1995.

\bibitem{Nolen_Ryzhik_09}
James Nolen and Lenya Ryzhik.
\newblock Traveling waves in a one-dimensional heterogeneous medium.
\newblock {\em Ann. Inst. H. Poincar{\'{e}} Anal. Non Lin{\'{e}}aire},
  26(3):1021--1047, 2009.

\bibitem{North_Ovaskainen_2007}
Ace North and Otso Ovaskainen.
\newblock Interactions between dispersal, competition, and landscape
  heterogeneity.
\newblock {\em Oikos}, 116(7):1106--1119, 2007.

\bibitem{Okubo_1989}
Akira Okubo, Philip~K. Maini, Mark~H. Williamson, and James~D. Murray.
\newblock On the spatial spread of the grey squirrel in {B}ritain.
\newblock {\em Proceedings of the Royal Society of London B: Biological
  Sciences}, 238(1291):113--125, 1989.

\bibitem{Osnas_2015}
Erik~E. Osnas, Paul~J. Hurtado, and Andrew~P. Dobson.
\newblock Evolution of pathogen virulence across space during an epidemic.
\newblock {\em The American Naturalist}, 185(3):332--342, 2015.

\bibitem{Ovaskainen_2008}
Otso Ovaskainen.
\newblock Analytical and numerical tools for diffusion-based movement models.
\newblock {\em Theoretical Population Biology}, 73(2):198 -- 211, 2008.

\bibitem{Perthame_2015}
Beno{\^{i}}t Perthame.
\newblock {\em Parabolic equations in biology}.
\newblock Lecture Notes on Mathematical Modelling in the Life Sciences.
  Springer, Cham, 2015.
\newblock Growth, reaction, movement and diffusion.

\bibitem{Phillips2019}
Ben~L. Phillips and T.~Alex Perkins.
\newblock Spatial sorting as the spatial analogue of natural selection.
\newblock {\em Theoretical Ecology}, Mar 2019.

\bibitem{Potts_Petrovskii}
Jonathan~R. Potts and Sergei~V. Petrovskii.
\newblock Fortune favours the brave: Movement responses shape demographic
  dynamics in strongly competing populations.
\newblock {\em Journal of Theoretical Biology}, 420:190--199, 2017.

\bibitem{Risler_2017}
Emmanuel Risler.
\newblock Competition between stable equilibria in reaction-diffusion systems:
  the influence of mobility on dominance.
\newblock {\em ArXiv e-prints}, mar 2017.

\bibitem{Rodrigo_Mimura_2001}
Marianito Rodrigo and Masayasu Mimura.
\newblock Exact solutions of reaction-diffusion systems and nonlinear wave
  equations.
\newblock {\em Japan J. Indust. Appl. Math.}, 18(3):657--696, 2001.

\bibitem{Schoener_1976}
Thomas~W. Schoener.
\newblock Alternatives to {L}otka-{V}olterra competition: models of
  intermediate complexity.
\newblock {\em Theoret. Population Biology}, 10(3):309--333, 1976.

\bibitem{Kawasaki_Shige}
Nanako Shigesada and Kohkichi Kawasaki.
\newblock {\em Biological invasions: theory and practice}.
\newblock Oxford University Press, UK, 1997.

\bibitem{SKT_79}
Nanako Shigesada, Kohkichi Kawasaki, and Ei~Teramoto.
\newblock Spatial segregation of interacting species.
\newblock {\em Journal of Theoretical Biology}, 79(1):83--99, 1979.

\bibitem{Skellam_1951}
John~Gordon Skellam.
\newblock Random dispersal in theoretical populations.
\newblock {\em Biometrika}, pages 196--218, 1951.

\bibitem{Soave_Zilio_2015}
Nicola Soave and Alessandro Zilio.
\newblock Uniform bounds for strongly competing systems: the optimal
  {L}ipschitz case.
\newblock {\em Arch. Ration. Mech. Anal.}, 218(2):647--697, 2015.

\bibitem{Van_Vuuren_Norbury_2000}
Jan~H. van Vuuren and John Norbury.
\newblock Conditions for permanence in well-known biological competition
  models.
\newblock {\em ANZIAM J.}, 42(2):195--223, 2000.

\bibitem{V_V_V}
Aizik~I. Volpert, Vitaly~A. Volpert, and Vladimir~A. Volpert.
\newblock {\em Traveling wave solutions of parabolic systems}, volume 140 of
  {\em Translations of Mathematical Monographs}.
\newblock American Mathematical Society, Providence, RI, 1994.
\newblock Translated from the Russian manuscript by James F. Heyda.

\bibitem{Volterra_1926}
Vito Volterra.
\newblock Variazioni e fluttuazioni del numero d{'}individui in specie animali
  conviventi.
\newblock {\em Mem. Accad. Naz. Lincei}, 2:31--113, 1926.

\bibitem{Xin_1993}
Jack~X. Xin.
\newblock Existence and nonexistence of traveling waves and reaction-diffusion
  front propagation in periodic media.
\newblock {\em J. Statist. Phys.}, 73(5-6):893--926, 1993.

\bibitem{Yatat_Couteront_Dumont}
V.~Yatat, P.~Couteron, and Y.~Dumont.
\newblock Spatially explicit modelling of tree–grass interactions in
  fire-prone savannas: A partial differential equations framework.
\newblock {\em Ecological Complexity}, 36:290 -- 313, 2018.

\bibitem{Yu_Zhao_2015}
Xiao Yu and Xiao-Qiang Zhao.
\newblock Propagation phenomena for a reaction{--}advection{--}diffusion
  competition model in a periodic habitat.
\newblock {\em Journal of Dynamics and Differential Equations}, pages 1--26,
  2015.

\bibitem{Zlatos_2017}
Andrej Zlato\v{s}.
\newblock Existence and non-existence of transition fronts for bistable and
  ignition reactions.
\newblock {\em Ann. Inst. H. Poincar{\'{e}} Anal. Non Lin{\'{e}}aire},
  34(7):1687--1705, 2017.

\end{thebibliography}
